\documentclass[10pt,a4paper]{amsart}
\usepackage{amssymb,graphicx,latexsym,diagrams}
\theoremstyle{plain}
\newtheorem{lemma}{Lemma}
\newtheorem{proposition}{Proposition}
\newtheorem{theorem}{Theorem}
\newtheorem{corollary}{Corollary}
\theoremstyle{definition}
\newtheorem{definition}{Definition}
\theoremstyle{remark}
\newtheorem{remark}{Remark}

\newtheorem{example}{Example}
\newarrowhead{+>} \rMapsto\leftmapsto\downmapsto\upmapsto
\newarrow{To} ---->
\newarrow{Mapsto} |--->
\newarrow{Noto} -----
\newarrow{Into} >--->
\newarrow{Onto} ----{>>}
\newarrow{Iso} >---{>>}
\newarrow{Inclus} C---> 
\newarrow{Equal} =====
\newarrow{Dash} {dash}{dash}{dash}{dash}{dash}
\newarrow{Dashto} {dash}{dash}{dash}{dash}{>}
\begin{document}
\title{Spin Borromean Surgeries}
\date{October 20, 2000 and, in revised form, March 26, 2002.}
\author[G. Massuyeau]{Gw\'ena\"el Massuyeau}
\address{Laboratoire de Math\'ematiques de Nantes\\
         UMR 6629 CNRS/Universit\'e de Nantes \\
         2, rue de la Houssini\`ere \\
         BP 92208 \\
         44322 Nantes Cedex 03\\ 
         France}
\email{massuyea@math.univ-nantes.fr}
\begin{abstract}
In 1986, Matveev defined  the notion of Borromean surgery 
for closed oriented $3$-manifolds and
showed that the equivalence relation generated by this move is characterized
by the pair (first betti number, linking form up to isomorphism).\\
We explain how this extends for $3$-manifolds 
with spin structure if we replace the
linking form by the quadratic form defined by the spin structure. 
We then show that the equivalence relation among closed spin $3$-manifolds 
generated by spin Borromean surgeries is characterized by the triple 
(first betti number, linking form up to isomorphism, Rochlin invariant
modulo~ $8$).
\end{abstract}
\maketitle
\section*{Introduction}
\label{sec:intro}
The notion of Borromean surgery was introduced by Matveev in \cite{Mat}
as an example of what he called a $\mathcal{V}$-surgery. Since then, 
this transformation has become the elementary move of Goussarov-Habiro finite
type invariants theory for oriented $3$-manifolds 
(\cite{Hab}, \cite{Gou}, \cite{GGP}).
Matveev showed that the equivalence relation, among closed oriented 
$3$-manifolds, generated by Borromean surgery 
is characterized by the pair:
\begin{displaymath}
(\beta_{1}(M),\textrm{ isomorphism class of }\lambda_{M}) ,
\end{displaymath}
where $\beta_{1}(M)$ is the first Betti number
of a $3$-manifold $M$ and 
\begin{displaymath}
TH_{1}(M;\mathbf{Z})  \otimes TH_{1}(M;\mathbf{Z})
\rTo^{\lambda_{M}} \mathbf{Q}/\mathbf{Z} 
\end{displaymath}
is its torsion linking form. This result gives a characterization
of degree $0$ invariants in Goussarov-Habiro theory 
for closed oriented $3$-manifolds.

As mentioned by Habiro and Goussarov, their finite type invariants theory 
(in short: ``FTI theory'') makes sense also for 
$3$-manifolds with spin structure because Borromean surgeries work well 
with spin structures (see \S\ref{sec:spinmove}). So, the question is:
\emph{what is the ``spin'' analogue of Matveev's theorem?}\\
For each closed spin $3$-manifold $(M,\sigma)$, a quadratic form 
\begin{displaymath}
TH_{1}(M;\mathbf{Z}) \rTo^{\phi_{M,\sigma}} \mathbf{Q}/\mathbf{Z}
\end{displaymath}
can be defined by many ways (see \cite{LL}, \cite{MS}, and
also \cite{Tur}, \cite{Gil}). 
The bilinear form associated to $\phi_{M,\sigma}$ is $\lambda_{M}$.  
Its Gauss-Brown invariant is equal to $-R_{M}(\sigma)$ modulo $8$, where
\begin{displaymath}
Spin(M) \rTo^{R_{M}} \mathbf{Z}_{16}
\end{displaymath}
is the Rochlin function of $M$, sending a spin structure $\sigma$ of $M$
to the modulo $16$ signature of a spin $4$-manifold which spin-bounds
$(M,\sigma)$. 
The main result of this paper is the following refinement of Matveev's theorem:
\begin{theorem}
\label{th:main}
Let $(M,\sigma)$ and $(M',\sigma')$ be connected closed spin $3$-manifolds.\\
Then, the following assertions are equivalent:
\begin{enumerate}
\item $(M,\sigma)$ and $(M',\sigma')$  can be obtained from one another
by spin Borromean surgeries, 
\item there exists a homology isomorphism
$f: H_{1}(M;\mathbf{Z}) \rTo H_{1}(M';\mathbf{Z})$ such that:
\begin{displaymath} 
\phi_{M,\sigma}=\phi_{M',\sigma'} \circ f|_{TH_{1}(M;\mathbf{Z})},
\end{displaymath}
\item $R_{M}(\sigma)=R_{M'}(\sigma')$ modulo $8$ and
there exists a homology isomorphism
$f: H_{1}(M;\mathbf{Z}) \rTo H_{1}(M';\mathbf{Z})$ such that:
\begin{displaymath} 
\lambda_{M}=\lambda_{M'} \circ f|_{TH_{1}(M;\mathbf{Z})}.
\end{displaymath}
\end{enumerate}
\end{theorem}
\noindent
The equivalence between assertions 2 and 3 will be the topological statement
of an algebraic fact:
a nondegenerate quadratic form on a finite Abelian group is determined,
up to isomorphism, by its associated bilinear form and its Gauss-Brown
invariant.\\

In \S\ref{sec:defmove} we recall Matveev's notion of $\mathcal{V}$-surgery.
With this background, we then recall the definition of Borromean surgery, and
give equivalent descriptions of other authors.\\ 
In \S\ref{sec:spinmove}, we clarify how all of these notions have 
to be understood in the spin case: in particular,
spin Borromean surgeries are introduced. 
As a motivation to Theorem \ref{th:main},  FTI for spin $3$-manifolds,
in the sense of Habiro and Goussarov, are then defined: 
the Rochlin invariant is shown
to be a finite type degree $1$ invariant.
It should be mentioned  that Cochran
and Melvin have proposed a different FTI theory in \cite{CM}, 
and have also refined
their theory to the case of spin manifolds.\\
\S\ref{sec:algquad} is of an algebraic nature. We recall some definitions
and results about quadratic forms on finite Abelian groups.
We also prove the above mentioned algebraic fact: the proof makes use of 
Kawauchi-Kojima classification of linking pairings.\\
\S\ref{sec:topquad} is the topological cousin of the former: we review 
the quadratic form $\phi_{M,\sigma}$. Starting from Turaev $4$-dimensional
definition in \cite{Tur}, we then give an intrinsic 
definition for $\phi_{M,\sigma}$  (no reference to dimension $4$).\\ 
\S\ref{sec:proof} is devoted to the proof of Theorem \ref{th:main}. It
goes as a refinement of the original proof by Matveev for the ``unspun'' case.
Last section will give some of its applications.
\subsection*{Acknowledgments}
The author wants to thank his advisor Pr. Christian Blanchet
who supervised this work, and Florian Deloup for conversations 
regarding quadratic forms.
He is also indebted to F\'elicie Pastore for correcting his written english.\\
\section{Borromean surgeries and equivalent moves}
\label{sec:defmove}
First of all, we want to recall the \emph{unifying} idea of $\mathcal{V}$-surgery
by Matveev in \cite{Mat}. This will allow us to have a more conceptual view 
of Borromean surgeries.
\subsection{Review of Matveev $\mathcal{V}$-surgeries}
\label{subsec:Matveev}
\begin{figure}[h]
\begin{center}
\includegraphics[width=5cm,height=4cm]{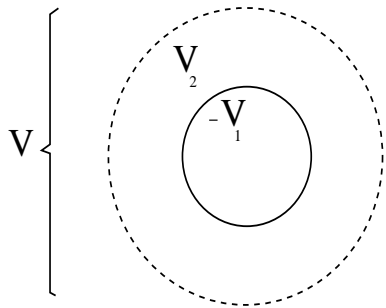}
\caption{What will be removed and what will be glued during the surgery.}
\label{fig:triple}
\end{center}
\end{figure}

We begin with some general definitions.
\begin{definition}
A \emph{Matveev triple} is a triple of oriented $3$-manifolds:
\begin{equation}
\label{eq:triple}
\mathcal{V}=(V,V_{1},V_{2}),
\end{equation}
where $V$ is closed and is the 
union of $-V_{1}$ and $V_{2}$ along their common boundary
$\partial V_{1}=\partial V_{2}$, as depicted in Figure~\ref{fig:triple}.\\
The triple $(-V,V_{2},V_{1})$ is called the \emph{inverse} of $\mathcal{V}$
and is denoted by $\mathcal{V}^{-1}$. 
\end{definition} 
Let now $M$ be a closed oriented $3$-manifold and let
$j: V_{1} \rTo M$ be an orientation-preserving
embedding. Form the following closed oriented $3$-manifold:
\begin{equation}
\label{eq:eqgensurg}
M'=M \setminus int\left(j(V_{1})\right) \cup_{j|_{\partial}} V_{2} .
\end{equation}
\begin{definition}
With the above notations, $M'$ is said to be obtained from $M$ by 
\emph{$\mathcal{V}$-surgery along $j$}.
\end{definition}
\noindent
Note that if $k$ denotes the embedding of $V_{2}$ in $M'$,
then $M$ is obtained from $M'$ by $\mathcal{V}^{-1}$-surgery along $k$.
\begin{definition}
Two Matveev triples $\mathcal{V}$ and $\mathcal{V'}$  are
said to be \emph{equivalent} if there exists an 
orientation-preserving diffeomorphism from
$V$ to $V'$ sending $-V_{1}$ to $-V'_{1}$, and $V_{2}$ to $V'_{2}$.
\end{definition}
\noindent
Note that, if the triples $\mathcal{V}$ and $\mathcal{V'}$ are equivalent,
then they have the same surgery effect.
\begin{example}  
\label{ex:triple}
Let $\Sigma_{g}$ denote the genus $g$ closed oriented
surface and let $H_{g}$ be the genus $g$ oriented handlebody. Then, 
each orientation-preserving diffeomorphism 
$f:\Sigma_{g} \rTo\Sigma_{g}$ leads to
a triple: 
\begin{displaymath}
\mathcal{V}_{f}:=\left((-H_{g})\cup_{f} H_{g},H_{g},H_{g}\right).
\end{displaymath}
A $\mathcal{V}_{f}$-surgery amounts to ``twist'' an embedded genus $g$
handlebody by $f$. For instance, from the standard genus one 
Heegaard decomposition of $\mathbf{S}^{3}$,
integral Dehn surgery is recovered.
\end{example}
\subsection{Review of Borromean surgeries}
\label{subsec:review_booromean}
\begin{figure}[h]
\begin{center}
\includegraphics[width=10cm,height=5cm]{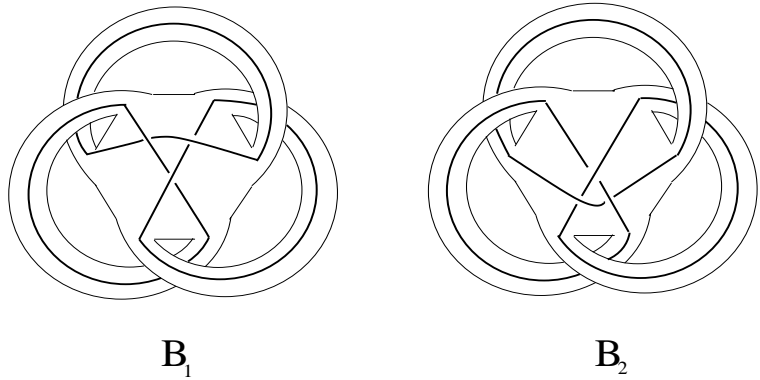}
\caption{The triple $\mathcal{B}$.}
\label{fig:halves}
\end{center}
\end{figure}

\begin{figure}[h]
\begin{center}
\includegraphics[height=3.5cm,width=4cm]{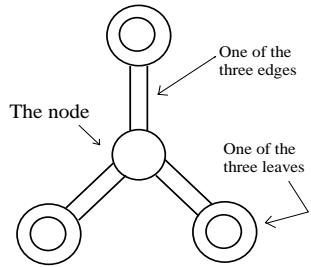} 
\caption{A Y-graph.}
\label{fig:Y}
\end{center}
\end{figure}

\begin{figure}[h]
\begin{center}
\includegraphics[height=4cm,width=5.5cm]{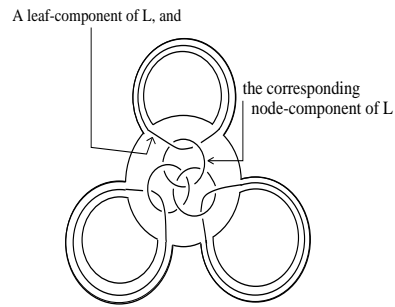} 
\caption{The surgery meaning of a $Y$-graph.}
\label{fig:L}
\end{center}
\end{figure}

The original Matveev's point of view was: 
\begin{definition}
A \emph{Borromean surgery} is a $\mathcal{B}$-surgery with: 
\begin{displaymath}
\mathcal{B}=
\left(B:=(-B_{1})\cup B_{2},B_{1},B_{2}\right),
\end{displaymath}
where the ``halves'' $B_{1}$ and $B_{2}$ are obtained from the genus $3$
handlebody by surgery  on three-component framed links as shown 
in Figure\footnote{Blackboard framing convention is used.} \ref{fig:halves}.
\end{definition}
We now recall Goussarov's notion of \emph{$Y$-surgery} in \cite{Gou}.
This move is equivalent to the \emph{$A_{1}$-move} of Habiro in \cite{Hab}.
\begin{definition}
\label{def:Ygraph}
A \emph{$Y$-graph} $G$ in a closed oriented $3$-manifold $M$ is an
(unoriented) embedding 
of the surface drawn in Figure \ref{fig:Y}, together with its decomposition
between  \emph{leaves}, \emph{edges} and \emph{node}.\\
The closed oriented $3$-manifold \emph{obtained 
from $M$ by $Y$-surgery along $G$} is: 
\begin{displaymath}
M_{G}:= (M \setminus int(N(G))) \cup (H_{3})_{L} ,
\end{displaymath}
where $N(G)\cong_{+} H_{3}$ is a regular neighbourhood of $G$ in $M$,
and $(H_{3})_{L}$ is the surgered handlebody
on the six-component link $L$ drawn on Figure \ref{fig:L}.\\
We call \emph{$Y$-equivalence} the equivalence relation among closed oriented 
$3$-manifolds generated by orientation-preserving 
diffeomorphisms and $Y$-surgeries. 
\end{definition} 
Note that a $Y$-surgery is a $\mathcal{Y}$-surgery if we call 
$\mathcal{Y}$ the triple:
\begin{displaymath}
\mathcal{Y}=\left(Y:=(-H_{3})\cup (H_{3})_{L},
 H_{3},(H_{3})_{L}\right), 
\end{displaymath}
the corresponding $Y$-graph gives the place
where the $\mathcal{Y}$-surgery must be performed.
\begin{lemma}
\label{lem:eqBY}
The Matveev triples $\mathcal{B}$ and $\mathcal{Y}$ are equivalent. Thus,
Borromean surgery is equivalent to $Y$-surgery.
\end{lemma}
\begin{proof} We will show that both of the triples 
$\mathcal{B}$ and $\mathcal{Y}$ are equivalent to a triple 
$\mathcal{V}_{h}$, defined by an orientation-preserving diffeomorphism
$h:\Sigma_{3} \rTo \Sigma_{3}$.\\
We start with the ``half'' $(H_{3})_{L}$ of Figure \ref{fig:L}:
handle-sliding of each node-component
over the corresponding leaf-component, followed by some isotopies
of framed links gives Figure \ref{fig:Linter},
where only part of the link is drawn. 
\begin{figure}[h]
\begin{center}
\includegraphics[height=4.5cm,width=5cm]{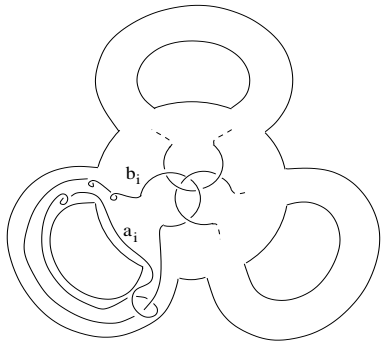}
\caption{The link $L$ of Figure \ref{fig:L} after some Kirby moves.}
\label{fig:Linter}
\end{center}
\end{figure}
Up to a $(+1)$-framing correction, the three  depicted components
$a_{i}$ can be normally pushed off at once towards the boundary: 
we obtain three disjoint curves $\alpha_{i}$ on $\partial H_{3}$.
Note that during this push-off, none of the three components $b_{i}$ 
is intersected. Then, the components $b_{i}$ can also be pushed off 
so that  the framing
correction is now $-1$: the result is a family of three disjoint curves 
$\beta_{i}$. After a convenient isotopy of the handles,
the curves $\alpha_{i}$ and $\beta_{i}$ 
can be  depicted as on Figure \ref{fig:curves}.
\begin{figure}[h]
\begin{center}
\includegraphics[height=5cm,width=10cm]{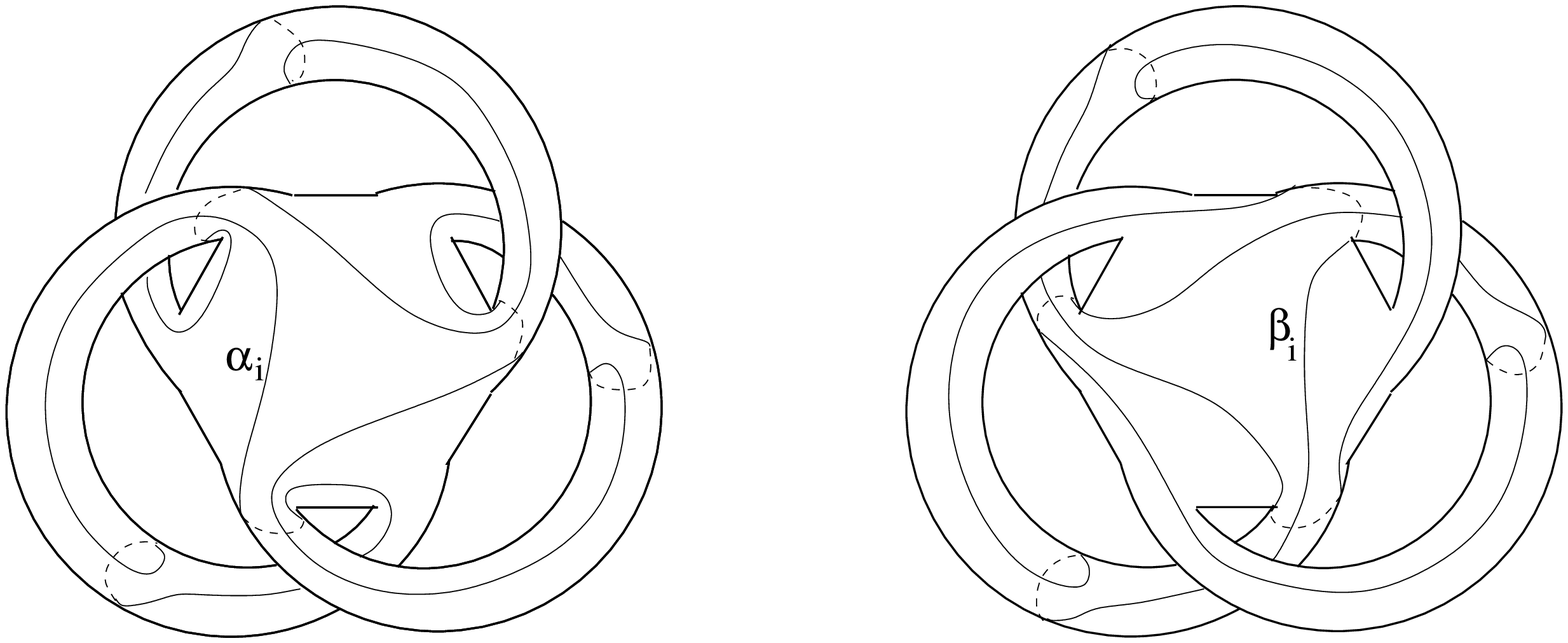}
\caption{The curves $\alpha_{i}$ (on the lefthand side) and $\beta_{i}$ 
(on the righthand side) resulting from
the push-off of the corresponding knots $a_{i}$ and $b_{i}$ of Figure 
\ref{fig:Linter}.}
\label{fig:curves}
\end{center}
\end{figure}

We define $h:=h_{b}^{-1} \circ h_{a}$, where $h_{a}$ and $h_{b}$ are
the following composites  of (commuting) Dehn twists:
\begin{displaymath}
h_{a}=\prod_{i=1}^{3} \tau_{\alpha_{i}} \quad \textrm{and} \quad
h_{b}=\prod_{i=1}^{3} \tau_{\beta_{i}}.   
\end{displaymath}
According to the Lickorish trick \cite[proof of Theorem 2]{Lic},  
a $\mathcal{Y}$-surgery is therefore
equivalent to a $\mathcal{V}_{h}$-surgery.

On the other hand, from Figure \ref{fig:halves}, 
we deduce that a $\mathcal{B}$-surgery 
is equivalent to a $\mathcal{V}_{h'}$-surgery where 
$h':=h_{b}\circ h_{a}^{-1}$. Let $m_{i}$ denote the meridian of the $i^{e}$
handle of $H_{3}$, for $i=1,2,3$. Then, the equation:
\begin{equation}
k^{-1}=k \circ \prod_{i=1}^{3} \tau_{m_{i}}^{2}
\end{equation}
holds for both $h_{a}$ and $h_{b}$ so that $h=h'$. 
\end{proof}
\begin{remark}[Fundamental remark]
\label{rem:fundamental}
Note that, in the proof of Lemma \ref{lem:eqBY}, the curve $\alpha_{i}$ is
homologous in the surface $\Sigma_{3}$ to the corresponding curve $\beta_{i}$
(look at Figure \ref{fig:curves}).
As a consequence, the diffeomorphism $h:\Sigma_{3} \rTo \Sigma_{3}$
belongs to the Torelli group. In the sequel, we will call $h$ the
\emph{Borromean diffeomorphism}.
\end{remark}
In order to have a complete understanding of these equivalent triples
$\mathcal{B}$, $\mathcal{Y}$ or $\mathcal{V}_{h}$, it remains to recognize
their underlying closed $3$-manifolds.
\begin{lemma}
\label{lem:whoisY}
The closed $3$-manifolds $B$ and $Y$, respectively defined by the triples 
$\mathcal{B}$ and $\mathcal{Y}$, are both homeomorphic to the $3$-torus 
$\mathbf{S}^{1}\times \mathbf{S}^{1}\times \mathbf{S}^{1}$.
\end{lemma}
\begin{proof} According to Lemma \ref{lem:eqBY},
$B$ and $Y$ are diffeomorphic. Let us
identify $Y$. Recall that $Y$ was defined as:
\begin{equation}
\label{eq:3manifoldY}
Y=(-H_{3})\cup (H_{3})_{L}.
\end{equation}
Write $L$ as $L_{a}\ \dot{\cup}\ L_{b}$ where $L_{a}$ (resp. $L_{b}$) is the
sublink containing the leaf(resp. node)-components.
Note that $L_{b}$ can be isotoped in $Y\setminus L_{a}$ 
to some Borromean rings contained in a $3$-ball disjoint from $L_{a}$: 
make $L_{b}$ leave the 
handlebody where it was lying, towards the handlebody with the minus sign
in (\ref{eq:3manifoldY}).
So $Y$ is obtained from $(-H_{3})\cup (H_{3})_{L_{a}}$ by surgery
on some Borromean rings contained in a little $3$-ball. The lemma
then follows from the fact that the latter 
is nothing but $\mathbf{S}^{3}$. 
\end{proof}
\section{Spin Borromean surgeries}
\label{sec:spinmove}
We now go into the world of spin $3$-manifolds. We refer
to \cite[Chapter IV]{Kir}  for an introduction to spin structures.
As a warming up, we recall a few facts in the next subsection.
\subsection{Glueing of spin structures}
\label{subsec:glueing}
Let $n\geq 2$, and let $M$ be a compact smooth oriented $n$-manifold
endowed with a Riemannian metric. Its bundle
of oriented orthonormal frames will be denoted by $FM$: it is a principal
$SO(n)$-bundle with total space  $E(FM)$ 
and with projection $p$.\\
Recall that if $M$ is spinnable, $Spin(M)$ can be thought of as:
\begin{displaymath}
Spin(M)=\left\{\sigma \in H^{1}\left(E(FM);\mathbf{Z}_{2}\right) \ : \ 
\sigma|_{\textrm{fiber}}\neq 0\in H^{1}\left(SO(n);\mathbf{Z}_{2}\right) \right\},
\end{displaymath}
and is essentially independant of the metric.
The set $Spin(M)$ is then an affine space over $H^{1}(M;\mathbf{Z}_{2})$, the 
action being defined by:
\begin{displaymath}
\forall x \in H^{1}(M;\mathbf{Z}_{2}),\ \forall \sigma \in Spin(M), \
x\cdot\sigma:= \sigma + p^{*}(x).
\end{displaymath}
\begin{lemma}
\label{lem:glueing}
For $i=1,2$, let $M_{i}$ be a compact smooth oriented $n$-manifold and let
$S_{i}$ be a submanifold of $\partial M_{i}$ with orientation
induced by $M_{i}$. 
Let also $f:S_{2} \rTo S_{1}$ be an orientation-reversing 
diffeomorphism and let $M=M_{1} \cup_{f} M_{2}$.\\
Assume that $M_{1}$ and $M_{2}$ are spinnable, that $S_{2}$ is
connected, and that the set: 
\begin{displaymath}
J=\left\{ (\sigma_{1},\sigma_{2})\in Spin(M_{1})\times Spin(M_{2}) \ : \
f^{*}\left(-\sigma_{1}|_{S_{1}}\right)=\sigma_{2}|_{S_{2}} \right\}
\end{displaymath}
is not empty. Then, $M$ is spinnable and the restriction map:
\begin{displaymath}
\left\{\begin{array}{rrcl}
Spin(M) & \rTo^r & Spin(M_{1})\times Spin(M_{2})\\
\sigma & \rMapsto & \left(\sigma|_{M_{1}},\sigma|_{M_{2}}\right),
\end{array}\right.
\end{displaymath}
is injective with $J$ as image.
\end{lemma}
\begin{proof}
For $i=1$ or $2$, let $F^+S_{i}$ be the principal $SO(n)$-bundle 
derived from $FS_{i}$ and the inclusion of groups:
\begin{diagram}
SO(n-1)& \rInclus^{g} &  SO(n).
\end{diagram}
We still denote by $g$ the 
canonical  map from $FS_{i}$ to $F^+S_{i}$. Then,
\begin{diagram}
H^{1}\left(E\left(F^+S_{i}\right);\mathbf{Z}_{2}\right)
& \rTo^{g^*} & H^{1}\left(E\left(FS_{i}\right);\mathbf{Z}_{2}\right)
\end{diagram}
is an isomorphism.
The bundle $F^+S_{i}$ can be identified with $FM_{i}|_{S_{i}}$.
In particular, there is an inclusion map 
$k_{i}:E\left(F^+S_{i}\right)\rInclus E(FM_i)$.
 
Moreover, the diffeomorphism $f$ induces a further identification:
\begin{diagram}
E\left(F^+S_{2}\right) & \rTo_{\simeq}^{f} & E\left(F^+S_{1}\right),
\end{diagram}
such that the total space $E(FM)$ is homeomorphic to the glueing:
\begin{displaymath}
E(FM_{1})\cup_{k_{1} f k_{2}^{-1}} E(FM_{2}).
\end{displaymath} 
We denote by $j_{i}$ the corresponding
inclusion of $E(FM_{i})$ into $E(FM)$. But now, by the Mayer-Vietoris
sequence, we have:
\begin{diagram}
H^{1}\left(E(FM);\mathbf{Z}_{2}\right) & 
\rTo^{\left(j_{1}^{*},j_{2}^{*}\right)} &
H^{1}\left(E(FM_1);\mathbf{Z}_{2}\right) \oplus  
H^{1}\left(E(FM_{2});\mathbf{Z}_{2}\right) \\
&&\dTo>{(k_{1}\circ f)^{*}-k_{2}^*}\\
&&H^{1}\left(E\left(F^+S_{2}\right);\mathbf{Z}_{2}\right). 
\end{diagram}
Note that $g^*\circ (k_{1}\circ f)^{*}$ sends each 
$\sigma_{1}\in Spin(M_{1})$ to  $f^{*}(-\sigma_{1}|_{S_{1}})$,
 while $g^*\circ k_2^*$ sends each 
$\sigma_{2}\in Spin(M_{2})$ to $\sigma_{2}|_{S_{2}}$. 
Note also that, since $S_{2}$ is connected, the map $(j_{1}^{*},j_{2}^{*})$
is injective. The whole lemma then follows from
these two remarks and from the exactness of the Mayer-Vietoris sequence. 
\end{proof}
\begin{definition}
\label{def:glueing}
With the notations and hypothesis of Lemma \ref{lem:glueing}, for
each $(\sigma_{1},\sigma_{2})\in J$, the unique spin structure 
of $M$ sent by $r$ to $(\sigma_{1},\sigma_{2})$ is called the 
\emph{glueing} of the spin structures $\sigma_{1}$ and $\sigma_{2}$, and
is denoted by $\sigma_{1} \cup \sigma_{2}$.
\end{definition}
\subsection{Spin  $\mathcal{V}$-surgeries}
\label{subsec:spinMatveev}
In some cases, a Matveev $\mathcal{V}$-surgery, whose definition has been
recalled in \S\ref{subsec:Matveev}, makes sense for spin $3$-manifolds.
\begin{definition}
A Matveev triple $\mathcal{V}=(V,V_1,V_2)$ 
is said to be \emph{spin-admissible}, 
if $\partial V_{1}=\partial V_{2}$ is connected, and
if the maps: 
\begin{displaymath}
H^{1}(V;\mathbf{Z}_{2})\rTo H^{1}(V_{1};\mathbf{Z}_{2}) \quad
 \textrm{and} \quad
H^{1}(V;\mathbf{Z}_{2})\rTo H^{1}(V_{2};\mathbf{Z}_{2}),
\end{displaymath}
induced by inclusions, are isomorphisms.
\end{definition}
Suppose now that $\mathcal{V}$ is a spin-admissible triple. Note that the
restriction maps: 
\begin{diagram}
Spin(V) & \rTo^{r_{1}} & Spin(-V_{1}) & \ \textrm{and} \ &
Spin(V) & \rTo^{r_{2}} & Spin(V_{2})
\end{diagram}
are then  bijective.
Let also $M$ be a closed oriented $3$-manifold and let $j:V_{1}\rTo M$ be
an orientation-preserving embedding. As in \S\ref{subsec:Matveev}, we denote
by $M'$ the result of the $\mathcal{V}$-surgery along $j$, and we want 
to define a \emph{canonical bijection}:
\begin{diagram}
Spin(M) & \rTo^{\Theta_{j,\mathcal{V}}} & Spin(M').
\end{diagram}
First, the embedding $j$ allows us to define the following map:
\begin{displaymath}
\left\{\begin{array}{rrcl}
Spin(M) & \rTo^{\Upsilon_{j,\mathcal{V}}}& Spin(V) \\
\sigma & \rMapsto & r_{1}^{-1}(-j^{*}(\sigma)).
\end{array}\right.
\end{displaymath}
From formula (\ref{eq:eqgensurg}) and from Definition 
\ref{def:glueing}, we can define $\Theta_{j,\mathcal{V}}(\sigma)$ 
as the following glueing:
\begin{displaymath}
\Theta_{j,\mathcal{V}}(\sigma):= \sigma|_{M\setminus int\left(j(V_{1})\right)}
\ \cup\ r_{2}\left(\Upsilon_{j,\mathcal{V}}(\sigma)\right).
\end{displaymath}
The inverse of $\Theta_{j,\mathcal{V}}$ is $\Theta_{k,\mathcal{V}^{-1}}$, where
$k$ denotes the embedding of $V_{2}$ in $M'$ .
\begin{definition}
With the above notations,
the spin manifold $\left(M',\Theta_{j,\mathcal{V}}(\sigma)\right)$ is said
to be obtained from $(M,\sigma)$ by 
\emph{spin $\mathcal{V}$-surgery along $j$}.
\end{definition} 
\begin{example} 
\label{ex:spin_triple}
Let $f: \Sigma_{g}\rTo \Sigma_{g}$ 
be an orientation-preserving diffeomorphism. 
Denote by $K$ the lagrangian subspace
of $H_{1}(\Sigma_{g};\mathbf{Z}_{2})$ span by the meridians. 
Then, as can be easily verified,
the triple $\mathcal{V}_{f}$ of Example \ref{ex:triple}
is spin-admissible if and only if $f_{*}(K)=K$. 
For instance, this condition is
satisfied when $f$ belongs to the Torelli modulo $2$ group.
\end{example} 
\begin{lemma}
\label{lem:spin_triple}
In the particular case of Example \ref{ex:spin_triple}, 
that is when $\mathcal{V}=\mathcal{V}_{f}$ with $f_{*}(K)=K$,
then for each $\sigma \in Spin(M)$, $\Theta_{j,\mathcal{V}}(\sigma)$ 
is the unique spin structure of $M'$ extending 
$\sigma|_{M\setminus int\left(j(V_{1})\right)}$.
\end{lemma}
\begin{proof}
In that case, $V_{2}$ is a handlebody and so 
$H^1\left(M',M\setminus int\left(j(V_{1})\right);
\mathbf{Z}_{2}\right)$ is zero. Therefore,
the restriction map
$Spin(M')\rTo Spin(M\setminus int(j(V_{1})))$ is injective.
\end{proof}
\subsection{Definition of $Y^{s}$-surgeries}
\label{subsec:defYspin}
From Lemma \ref{lem:eqBY} and from Remark \ref{rem:fundamental} above,
we have learnt that both of the triples 
$\mathcal{B}$ and $\mathcal{Y}$ are equivalent to the triple 
$\mathcal{V}_{h}$ where $h$ is the Borromean diffeomorphism which belongs
to the Torelli group. So, by
Example \ref{ex:spin_triple}, they are spin-admissible and the following
definition makes sense: 
\begin{definition}
A \emph{$Y^{s}$-surgery}, or equivalently a \emph{spin Borromean surgery},
is the surgery move among closed spin $3$-manifolds defined 
equivalently by the triples $\mathcal{Y}$ or $\mathcal{B}$.
We call \emph{$Y^{s}$-equivalence} the equivalence relation
among them generated by spin diffeomorphisms
and $Y^{s}$-surgeries.
\end{definition} 
Let $(M,\sigma)$ be a closed spin $3$-manifold 
and let $G$ be a $Y$-graph in $M$. 
The $Y^{s}$-surgery along $G$ gives a new spin manifold which will be  
denoted by:
\begin{displaymath}
(M_{G},\sigma_{G}).
\end{displaymath}
Let $j:H_{3} \rTo N(G)$ be an embedding of the genus 3
handlebody onto a regular neighbourhood of $G$ in $M$. Then, 
\begin{equation}
\label{eq:how_to_think_M'}
M_{G}\cong  M\setminus int\left(N(G)\right) \cup_{j|_{\partial}\circ h} H_{3},
\end{equation}
and according to Lemma \ref{lem:spin_triple}, 
$\sigma_{G}$ is the unique spin structure of $M_{G}$ extending 
$\sigma|_{M\setminus int\left(N(G)\right)}$.
\subsection{Goussarov-Habiro FTI theory for spin $3$-manifolds}
\begin{lemma}
\label{lem:associativity}
Let $(M,\sigma)$ be a closed spin $3$-manifold and let $G$ and $H$
be disjoint $Y$-graphs in $M$. Then,
up to diffeomorphism of manifolds with spin structure, 
\begin{displaymath}
\left((M_{G})_{H},(\sigma_{G})_{H}\right)=
\left((M_{H})_{G},(\sigma_{H})_{G}\right).
\end{displaymath}
\end{lemma}
\begin{proof}
The equality $(M_{G})_{H}=(M_{H})_{G}$ is obvious. By construction, 
both of $(\sigma_{G})_{H}$ and $(\sigma_{H})_{G}$ are extensions of 
$\sigma|_{M\setminus \left(N(G)\cup N(H)\right)}$. The lemma then follows
from the fact that the restriction map:
\begin{diagram}
Spin\left((M_{G})_{H}\right) & \rTo & 
Spin\left(M\setminus \left(N(G)\cup N(H)\right)\right)
\end{diagram}
is injective since the relative cohomology group 
$H^{1}\left((M_{G})_{H},M\setminus \left(N(G)\cup N(H)\right);
\mathbf{Z}_{2}\right)$ is zero.
\end{proof}
Let $S=\{G_{1},\dots,G_{s}\}$ be a family of disjoint $Y$-graphs in
a closed $3$-manifold $M$ with spin structure $\sigma$. 
Lemma \ref{lem:associativity}
says that \emph{$Y^{s}$-surgery along the family} $S$ is well-defined. We denote
the result by $(M_{S},\sigma_{S})$.
Following Habiro and Goussarov definition of a finite type invariant 
(\cite{Hab}, \cite{Gou}), we can now define:
\begin{definition}
\label{def:spinFTI}
Let $A$ be an Abelian group and let $\lambda$ be an $A$-valued invariant
of $3$-manifolds with spin structure. Then, $\lambda$ is an 
\emph{invariant of degree at most $n$} if for any closed spin $3$-manifold 
$(M,\sigma)$ and any family $S$ of at least $n+1$ $Y$-graphs in $M$, 
the following identity holds:
\begin{equation}
\sum_{S'\subset S} \lambda(M_{S'},\sigma_{S'})=0 \in A,
\end{equation}
where the sum is taken over all subfamilies $S'$ of $S$.
Moreover, $\lambda$ is \emph{of degree $n$} if it is of degree
at most $n$, but is not of degree at most $n-1$.
\end{definition}
\begin{remark}
Note that the degree $0$ invariants are precisely those invariants which are
constant on each $Y^{s}$-equivalence class. So, the refined Matveev
theorem will quantify how powerful they can be.
\end{remark}
The next subsection will provide us some examples of invariants
which are  finite type in the sense of Definition \ref{def:spinFTI}.
\subsection{Rochlin invariant under $Y^{s}$-surgery}
\label{subsec:Rochlin}
\begin{proposition}
\label{prop:Rochlin}
Let $(M,\sigma)$ be a closed spin $3$-manifold, and let $G$ be a
$Y$-graph in $M$.
Then, the following formula holds:
\begin{displaymath}
R_{M_{G}}(\sigma_{G})=R_{M}(\sigma)
+ R_{\mathbf{S}^{1}\times \mathbf{S}^{1}\times\mathbf{S}^{1}}
\left(\Upsilon_{G}(\sigma)\right)
\in \mathbf{Z}_{16} ,
\end{displaymath}
where the map $\Upsilon_{G}: Spin(M) \rTo 
Spin(\mathbf{S}^{1}\times \mathbf{S}^{1}\times\mathbf{S}^{1})$, induced by 
the $\mathcal{Y}$-surgery along $G$,
has been defined in \S\ref{subsec:spinMatveev}.\\
\end{proposition}
\begin{proof}  
According to Lemma \ref{lem:whoisY}, we can think 
of the $3$-torus as: 
\begin{equation}
\label{eq:3torus}
\mathbf{S}^{1}\times \mathbf{S}^{1}\times\mathbf{S}^{1} 
=(-H_{3})\cup_{h} H_{3}.
\end{equation}
The surgered manifold $M_{G}$ will be thought of 
concretely as in (\ref{eq:how_to_think_M'}).\\
Pick  a spin $4$-manifold $W$ spin-bounded by $(M,\sigma)$, and a
spin $4$-manifold $H$ spin-bounded by the $3$-torus with 
$\Upsilon_{G}(\sigma)$ as a spin structure.
Glue the ``generalized'' handle $H$ to $W$ along 
the first handlebody of the $3$-torus in decomposition (\ref{eq:3torus}),
using $j$ as  glue.
We obtain a $4$-manifold $W'$. 
Orient $W'$ coherently with $H$ and $W$, and then 
give to $W'$ the spin structure obtained by glueing those of $H$
and $W$ (see Definition \ref{def:glueing}).
It follows from definitions that  the spin-boundary
of $W'$ is $(M_{G},\sigma_{G})$.\\
According to Wall theorem on non-additivity of the signature 
(see \cite{WalIII}), we have:
\begin{equation}
\label{eq:signatures}
sgn(W')=sgn(W)+sgn(H) - \textrm{ correcting term}.
\end{equation}
The involved correcting term is the signature of a real 
bilinear symmetric form explicitely described by Wall.
It is defined by means of the intersection form in $\Sigma_{3}$, 
with domain:
\begin{displaymath}
V=\frac{A\cap(B+C)}{A\cap B+A\cap C},
\end{displaymath}
where $A$,$B$,$C$ are subspaces of $H_{1}(\Sigma_{3};\mathbf{R})$ defined
to be respectively the kernels of:
\begin{diagram}
&H_{1}(\Sigma_{3};\mathbf{R}) & \rTo^{\left(j|_{\partial}\right)_{*}} &
H_{1}(\partial N(G);\mathbf{R}) & \rTo^{i_{*}} &
H_{1}(M\setminus int(N(G));\mathbf{R}),\\
&H_{1}(\Sigma_{3};\mathbf{R}) & &  \rTo^{i_{*}} & & H_{1}(H_{3};\mathbf{R}),\\
\textrm{and:} \quad &H_{1}(\Sigma_{3};\mathbf{R}) & \rTo^{(h^{-1})_{*}}& 
H_{1}(\Sigma_{3};\mathbf{R}) & \rTo^{i_{*}} & H_{1}(H_{3};\mathbf{R}).
\end{diagram}
No matter who is $A$, since the Borromean diffeomorphism $h$ 
lies in the Torelli group, we certainly
have $B=C$. The space $V$ then vanishes and so does the correcting term.
The announced equality then follows by taking equation
(\ref{eq:signatures}) modulo $16$. 
\end{proof}
\begin{corollary}
\label{cor:Rochlin}
The Rochlin invariant is a degree $1$ invariant of closed 
spin $3$-manifolds for Goussarov-Habiro theory, 
and its modulo $8$ reduction is of degree $0$. 
\end{corollary}
\begin{remark}
In Cochran-Melvin  theory, the Rochlin invariant is a degree $3$ finite type
invariant (see \cite[Prop. 6.2]{CM}).
\end{remark}
\begin{proof}[Proof of Corollary \ref{cor:Rochlin}]
Last statement is clear from Proposition \ref{prop:Rochlin} and
from the fact that the Rochlin function of the
$3$-torus takes values in $\{0,8\}\subset \mathbf{Z}_{16}$. Let us show that
the Rochlin invariant is at most of degree $1$. 
Take a closed spin $3$-manifold $(M,\sigma)$ and 
two disjoint $Y$-graphs $G$ and $H$ in $M$. According to
Proposition \ref{prop:Rochlin}, in order to show that:
\begin{equation}
R_{M}(\sigma)-R_{M_{G}}(\sigma_{G})-R_{M_{H}}(\sigma_{H})
+R_{M_{G,H}}(\sigma_{G,H})=0,
\end{equation}
it suffices to show that:
\begin{equation}
\Upsilon_{G}(\sigma)=\Upsilon_{G} (\sigma_{H})
\in Spin(\mathbf{S}^{1}\times \mathbf{S}^{1}\times\mathbf{S}^{1}),
\end{equation}
where the left $\Upsilon_{G}$ is defined by $G\subset M$ and the
right $\Upsilon_{G}$ is defined by $G\subset M_{H}$.
But this follows from definition of the maps $\Upsilon_{G}$ and from 
the fact that $\sigma_{H}$ extends
$\sigma|_{M\setminus N(H)}$.\\
It remains to show that the Rochlin invariant is not of degree $0$ 
(and so it will be 
``exactly'' of degree $1$). For instance, all of the spin structures of
the $3$-torus are related one to another by $Y^{s}$-surgeries 
(\emph{Cf} Example \ref{ex:torus} below)
and so are not distinguished one to another by degree $0$ invariants.
But Rochlin distinguishes one of them from the others. 
\end{proof}
\subsection{$Y^{s}$-surgeries through 
   surgery presentations on $\mathbf{S}^{3}$}
\label{subsec:surgery_presentations}
\subsubsection*{We first fix some notations.} 
We call $V_{L}$ the $3$-manifold obtained from $\mathbf{S}^{3}$
by surgery along a ordered oriented framed link 
$L=(L_{1},\dots,L_{l})$ of length $l$, and $W_{L}$ the 
corresponding $4$-manifold obtained from $\mathbf{B}^{4}$ by attaching
$2$-handles and sometimes called the \emph{trace} of the surgery.
Let also $B_{L}=(b_{ij})_{i,j=1,\dots,l}$ be the linking matrix of $L$.\\
Recall that $H_{2}(W_{L};\mathbf{Z})$ is free Abelian of rank $l$. For
each $i=1,\dots,l$, choose a Seifert surface of $L_{i}$ in $\mathbf{S}^{3}$
and push it off into the interior of $\mathbf{B}^{4}$: denote the result by 
$P_{i}$. Then, glue $P_{i}$ to the core of the $i^{e}$ $2$-handle to obtain 
a closed surface $S_{i}$. A basis of $H_{2}(W_{L};\mathbf{Z})$ is then
given by $([S_{1}],\dots,[S_{l}])$.
\subsubsection*{We now recall  a nice combinatorial way 
to describe spin structures of $V_{L}$.} This is by means of
the so-called ``characteristic solutions of $B_{L}$'' or, equivalently,
``characteristic sublinks of $L$''.\\
A vector $s\in (\mathbf{Z}_{2})^{l}$ or, equivalently, the sublink of $L$ 
containing the components $L_{i}$ such that $s_{i}=1$, are said to be 
\emph{characteristic} if the following equation is satisfied:
\begin{equation}
\label{eq:charsol}
\forall i=1,\dots,l \ ,  \quad \sum_{j=1}^{l}b_{ij}\cdot s_{j}= b_{ii}
\in \mathbf{Z}_{2}.
\end{equation}
We denote by $\mathcal{S}_{L}$ the subset of $(\mathbf{Z}_{2})^{l}$ comprising
the characteristic solutions of $B_{L}$. There is a bijection: 
\begin{diagram}
Spin(V_{L}) & \rTo^{\simeq} & \mathcal{S}_{L} 
\end{diagram}
which is defined by the following composition:
\begin{diagram}
Spin(V_{L}) & \rTo^{o} & H^{2}(W_{L},V_{L};
\mathbf{Z}_{2}) & \rTo^{P}_{\simeq} &
H_{2}(W_{L};\mathbf{Z}_{2}) & \rTo_{\simeq} & (\mathbf{Z}_{2})^{l}
\end{diagram}
where $o$ sends any $\sigma\in Spin(V_L)$ 
to the obstruction to extend $\sigma$ to the
whole of $W_{L}$, $P$ is the Poincar\'e duality isomorphism 
and the last map is
defined by the basis $([S_{1}],\dots,[S_{l}])$.
With this combinatorial description, 
Kirby theorem can be refined to closed spin
$3$-manifolds (see \cite{Bl}).\\ 
The following lemma, more general than needed, will allow  
us to enunciate in those terms the effect of a $Y^{s}$-surgery.
\begin{lemma}
\label{lem:spintech}
Let $L\cup K$ be the ordered union of two ordered oriented framed links in  
$\mathbf{S}^{3}$ and let $H\subset S^{3}$ be an embedded handlebody
such that $K$ is contained in the interior of $H$, $L$
is disjoint from $H$, and $H_{K}$ is a  $\mathbf{Z}_{2}$-homology handlebody.

Suppose that $\sigma \in Spin\left(V_{L \cup K}\right)$ is represented by
a characteristic solution $s\in (\mathbf{Z}_{2})^{l+k}$ of $B_{L\cup K}$
satisfying the following two properties:
\begin{enumerate}
\item $(s_{1},\dots,s_{l}) \in (\mathbf{Z}_{2})^{l}$ is
a characteristic solution of $B_L$,
\item for all $i \in \{1,\dots,k\}$ such that $s_{l+i} \neq 0$,
the component $K_{i}$ bounds a Seifert surface \emph{within} $H$. 
\end{enumerate}
Then, the restricted spin structure: 
\begin{displaymath}
\sigma_{|} \in Spin(V_{L \cup K} 
\setminus int(H_{K}))=Spin(V_{L} \setminus int(H)),
\end{displaymath}
extends to the spin structure of $V_{L}$
represented by $(s_{1},\dots,s_{l}) \in (\mathbf{Z}_{2})^{l}$.
\end{lemma}
\begin{proof} In the following, all (co)homology groups are assumed to be
with coefficients in $\mathbf{Z}_{2}$.
We use the above fixed notations.\\
Let us consider the map:
\begin{diagram}
Spin\left(V_{L} \setminus H\right) & \rTo^o & 
H^{2}\left(W_{L},V_{L}\setminus H\right),
\end{diagram}
where $o(\alpha)$ is the obstruction to extend any 
$\alpha \in Spin(V_L\setminus H)$ to the 
whole of $W_{L}$. Let also
$\delta^{*}:H^{1}\left(V_{L}\setminus H\right) \rTo 
H^{2}\left(W_{L},V_{L}\setminus H\right)$
denote the connecting homomorphism 
for the pair $\left(W_{L},V_{L}\setminus H\right)$.
Note that the following equation holds:
\begin{equation}
\forall x \in H^{1}\left(V_{L}\setminus H\right),
\forall \alpha \in Spin(V_L\setminus H),\
o(x\cdot\alpha)=o(\alpha)+\delta^{*}(x).
\end{equation}
Since $\delta^{*}$ is injective, it follows that $o$ is injective.\\
The same map $o$ can be defined for $V_{L}$ relatively to $W_{L}$,
and for $V_{L\cup K}$ and $V_{L\cup K} \setminus H_{K}$ relatively to
$W_{L \cup K}$.
We have thus the following commutative diagram:
\begin{diagram}
Spin(V_{L\cup K}) & \rInclus^o & H^{2}(W_{L\cup K},V_{L \cup K}) & 
\rTo^{\simeq}_P & H_{2}(W_{L\cup K}) & \rDash & \HmeetV\\
\dInclus &  & \dInclus & & \dInclus & & \dDash \\
Spin(V_{L\cup K} \setminus H_{K}) & \rInclus^o & 
 H^{2}(W_{L\cup K},V_{L \cup K}\setminus H_{K} ) & \rTo^{\simeq}_P &
 H_{2}(W_{L\cup K},H_{K}) & & r\\
\dEqual & & \dTo & & & & \dDash\\
Spin(V_{L}\setminus H) & \rInclus^o & H^{2}(W_{L},V_{L}\setminus H) &
\rTo^{\simeq}_P & H_{2}(W_{L},H) & \lDashto & \HmeetV\\
\uInclus & & \uInclus & & \uInclus>{i_{*}} \\
Spin(V_{L}) & \rInclus^o & H^{2}(W_{L},V_{L}) & 
\rTo^{\simeq}_P & H_{2}(W_{L}) \\
\end{diagram}  
where the letter $P$ stands for a 
Poincar\'e duality isomomorphism, the vertical
arrows are induced by inclusions and the map $r$ is defined by 
planar commutativity.\\
From intersection theory, we deduce that: 
\begin{equation}
\label{eq:calculdeh}
\left\{ \begin{array}{lr}
\forall i \in \{1,\dots,l\}, & r\left([S_{i}]\right)=[S_{i}], \\
\forall i \in \{l+1,\dots,l+k\},  & r\left([S_{i}]\right)=[P_{i}].
\end{array} \right.
\end{equation} 

Let now $\sigma$ be a spin structure of $V_{L\cup K}$ 
such that the corresponding characteristic solution $s$
of $B_{L\cup K}$ satisfies the conditions 1 and 2 of the lemma.\\ 
We define $\tilde{s} := (s_{1},\dots,s_{l}) \in (\mathbf{Z}_{2})^{l}$. 
By hypothesis 1, there exists a unique spin structure $\tilde{\sigma}$ 
of $V_{L}$ with $\tilde{s}$ as associated characteristic solution of $B_L$.
We want to show that 
$\sigma_{|}=\tilde{\sigma}_{|}$. Diagram chasing shows that proving
$r\circ P \circ o(\sigma)= i_{*}\circ P \circ o (\tilde{\sigma})$ should
suffice.
This follows from hypothesis 2, formulas (\ref{eq:calculdeh})
and from the fact that $P\circ o (\sigma)=s$ 
and $P\circ o (\tilde{\sigma})=\tilde{s}$. 
\end{proof}
\subsubsection*{Let us now come back to the case of a $Y^{s}$-surgery.}
Let $M$ be a closed oriented $3$-manifold and let $M_{G}$ be obtained from $M$
by surgery along a $Y$-graph $G$. According to \S\ref{subsec:spinMatveev},
$\mathcal{Y}$-surgery along $G$ induces a bijective map:
\begin{diagram}
Spin(M) & \rTo^{\Theta_{G}} & Spin(M_{G}).
\end{diagram}
With the notations of \S\ref{subsec:defYspin}, each $\sigma$ is
sent by $\Theta_{G}$ to $\sigma_{G}$.

\begin{figure}[h]
\begin{center}
\includegraphics[height=4cm,width=5cm]{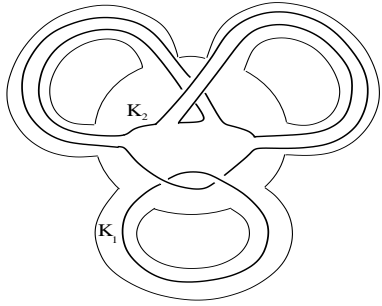}
\caption{$Y$-surgery as surgery along a $2$-component link.}
\label{fig:two-comp}
\end{center}
\end{figure}

Suppose now that we are given a surgery presentation $M=V_{L}$ of $M$. 
Isotope the graph $G$
in $M$ to make it disjoint from the dual of $L$, so that $G$ is in
$\mathbf{S}^{3}\setminus L$. Let $H$ be a regular neighbourhood of $G$.
A few Kirby calculi, inside $H$, show that surgery along this $Y$-graph 
is equivalent to surgery 
on the two-component link $K$ of Figure \ref{fig:two-comp}. We prefer this 
unsymmetric link to Figure \ref{fig:L} because 
of the fewer components.
\\
We then have $M_{G}=V_{L\cup K}$.
The linking matrix of $L\cup K$, 
when $K$ is appropriately oriented, looks like:
\begin{displaymath}
B_{L\cup K}=\left(
\begin{array}{ccc|cc}
& & & x_{1} & 0\\
& B_{L} & & \vdots & \vdots\\
& & & x_{l} & 0\\
\hline
x_{1} &\cdots & x_{l} & x & 1\\
0 & \cdots & 0 & 1 & 0 
\end{array} \right).
\end{displaymath}
Writing the characteristic condition (\ref{eq:charsol}), we find
that each characteristic solution of $B_{L\cup K}$ is of the form:
\begin{equation}
\label{eq:Yspin-chsol}
T_G(s):=\left(s_{1},\dots,s_{l},0,x+\sum_{i=1}^{l}x_{i} s_{i}\right)
\in (\mathbf{Z}_2)^{l+2},
\end{equation}
where $s=(s_{1},\dots,s_{l})$ must be characteristic for $L$.
Equation (\ref{eq:Yspin-chsol}) then defines a combinatorial bijection: 
\begin{diagram}
\mathcal{S}_{L} & \rTo^{T_{G}}_\simeq & \mathcal{S}_{L\cup K}.
\end{diagram}
\begin{lemma} 
With the above notations, the map $T_{G}$ is a combinatorial version
of the map $\Theta_{G}$ in terms of characteristic solutions for surgery 
presentations on $\mathbf{S}^{3}$. More precisely, the following
diagram is commutative: 
\begin{diagram}
\mathcal{S}_{L} & \rTo^{T_{G}} & \mathcal{S}_{L\cup K} \\
\uTo<{\simeq}& & \uTo>{\simeq} \\
Spin(M) & \rTo^{\Theta_{G}} & Spin(M_{G})
\end{diagram}
\end{lemma}
\begin{proof}
This follows from the definitions and from Lemma \ref{lem:spintech}:
note that $K_{2}$ is nul-homologous in $H$, and that here $H_{K}$ is merely
a handlebody.
\end{proof}
\begin{definition}
\label{def:simpleY}
Let $M=V_{L}$ be a surgery presentation of a $3$-manifold $M$ 
on $\mathbf{S}^{3}$, and let $G$ be a $Y$-graph in $M$.
Then, $G$ is said to be \emph{simple} (with respect to this surgery 
presentation), if $G$ can be isotoped in $M$ so that, in 
$\mathbf{S}^{3}$, its leaves bound disjoint discs, each intersecting $L$ 
in exactly one point.  
\end{definition} 
\begin{corollary}
\label{cor:simpleY}
For a $Y^{s}$-surgery along a simple $Y$-graph,
the spin-diffeomorphism of Figure \ref{fig:simpleY} holds.
\begin{figure}[!h]
\begin{center}
\includegraphics[width=10cm,height=4.5cm]{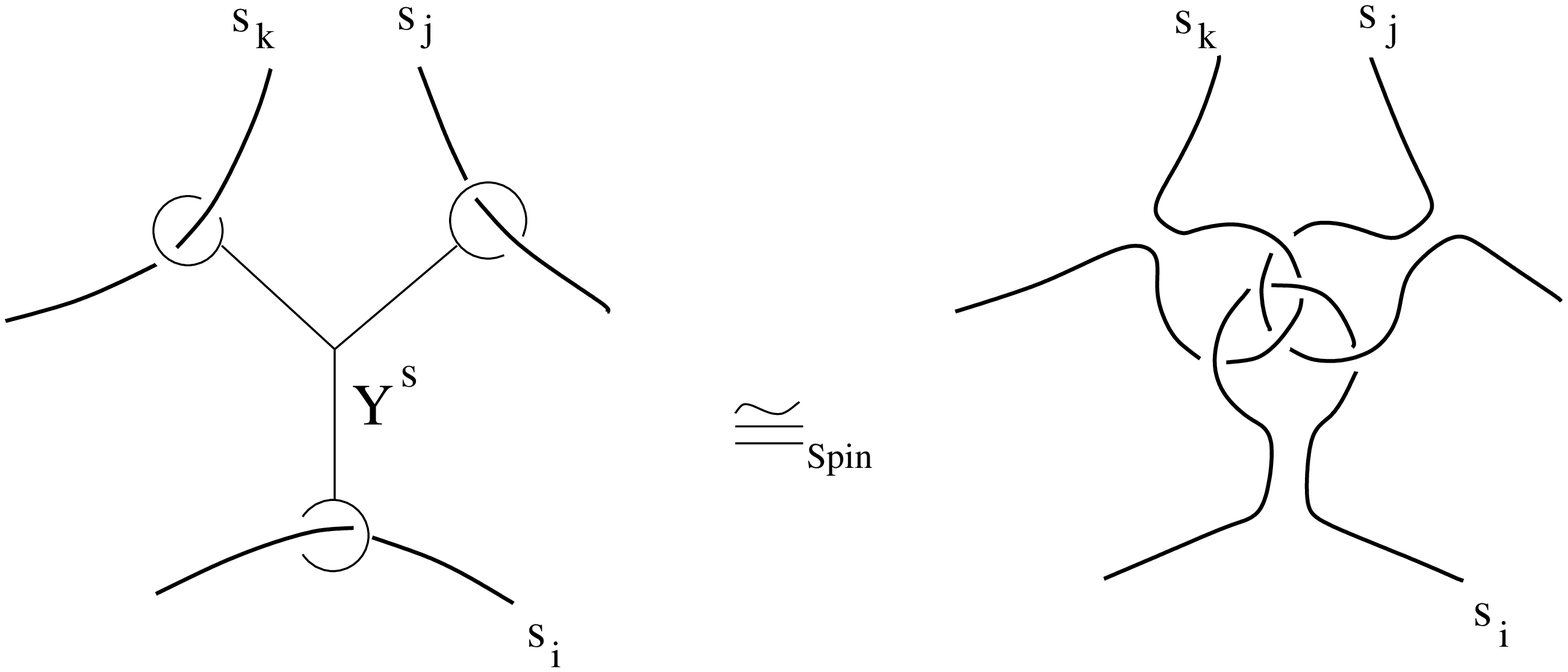}
\caption{A simple $Y^s$-surgery.}
\label{fig:simpleY}
\end{center}
\end{figure}
\end{corollary}
\begin{proof}  Replace in the lhs of Figure
\ref{fig:simpleY}, this simple $Y$-graph by
the 2-component link of Figure \ref {fig:two-comp} such that $K_{1}$
links the $i^{e}$ component of $L$ and use equation (\ref{eq:Yspin-chsol})
to obtain the intermediate link of Figure \ref{fig:simpleYinter}.
\begin{figure}[!h]
\begin{center}
\includegraphics[width=3cm,height=3cm]{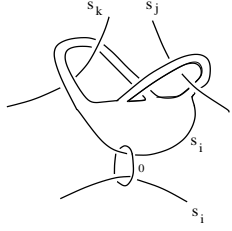}
 \caption{Some intermediate link and characterictic solution.}
\label{fig:simpleYinter}
\end{center}
\end{figure}
Perform then some spin Kirby moves to obtain the rhs
of Figure \ref{fig:simpleY}.
\end{proof}
\begin{example} 
\label{ex:torus}
As a consequence of
Corollary \ref{cor:simpleY}, the Lie spin structure of the $3$-torus
is $Y^{s}$-equivalent to the seven other ones.
\end{example} 
\begin{proof}  The two spin structures of $\mathbf{S}^{1} \times
\mathbf{S}^{2}$ are equivalent, so are the eight ones of
$\sharp^{3} \mathbf{S}^{1} \times \mathbf{S}^{2}$. Furthermore,
$\sharp^{3} \mathbf{S}^{1} \times \mathbf{S}^{2}$ can be obtained
from $\mathbf{S}^{3}$ by surgery along a trivial $0$-framed three-component
link. Surgery on the $0$-framed Borromean rings gives rise to 
the $3$-torus $\mathbf{S}^{1} \times \mathbf{S}^{1} \times \mathbf{S}^{1}$,
and this link can be obtained from the trivial 
link by a simple $Y$-surgery.
\end{proof}
\section{Quadratic forms on finite Abelian groups}
\label{sec:algquad}
\subsection{Linking pairings, quadratic forms
on finite Abelian groups and their presentations.}
\label{subsec:presentations}
We recall here standard algebraic constructions: notations are that
of Deloup in \cite{DelI}, where a brief review of the subject
can be found.

\begin{definition}
A \emph{linking pairing} on a finite Abelian group $G$ is a nondegenerate
symmetric bilinear map $b: G \times G \rTo \mathbf{Q}/\mathbf{Z}$.\\
A \emph{quadratic form} on $G$ is a map $q:G \rTo \mathbf{Q}/\mathbf{Z}$ 
such that the map $b_{q}:G \times G \rTo \mathbf{Q}/\mathbf{Z}$ defined by
$b_{q}(x,y)=q(x+y)-q(x)-q(y)$ is bilinear, and such that $q$ satisfies:
$\forall x \in G,\ q(-x)=q(x)$.
$q$ is said to be \emph{nondegenerate} when the associated 
bilinear form $b_{q}$ is a linking pairing.
\end{definition}
Let now $(F,f)$ be a symmetric bilinear form on a free finitely
generated Abelian group $F$. We denote by
$ad_{f}: F \rTo Hom(F,\mathbf{Z})$ the adjoint map, and by 
$a_{f}: F \otimes \mathbf{Q} \rTo 
Hom(F,\mathbf{Q})$ its rational extension. Form:
\begin{displaymath}
K_{f}:=\frac{Hom(F,\mathbf{Z}) \cap Im(a_{f})}{Im(ad_{f})}.
\end{displaymath}
Note that $K_{f}=T\left(Coker(ad_{f})\right)$, the torsion subgroup of
$Coker(ad_{f})$.
We now define a linking pairing: 
\begin{displaymath}
K_{f} \times K_{f} \rTo^{L_f} \mathbf{Q}/ \mathbf{Z}
\end{displaymath}
by the formula: 
\begin{equation}
L_{f}(\bar{x},\bar{y})=x_{\mathbf{Q}}(\tilde{y})\quad
 mod \quad 1,
\end{equation}
where $x,y \in Hom(F,\mathbf{Z})\cap Im(a_{f})$, 
$\tilde{y} \in F \otimes \mathbf{Q}$ is such that $a_{f}(\tilde{y})=y$,
and $x_{\mathbf{Q}} \in Hom(F\otimes \mathbf{Q}, \mathbf{Q})$ is the rational
extension of $x$.
$(F,f)$ is said to be a \emph{presentation} of the linking
pairing $(K_{f},L_{f})$.

Suppose now that the form $(F,f)$ comes equipped with a \emph{Wu class},
that is an element $w\in F$ such that $\forall x \in F,f(w,x)=f(x,x)
 \quad mod \quad 2$.
We can then define a quadratic form over $L_{f}$, denoted by:
\begin{displaymath}
K_{f} \rTo^{\phi_{f,w}} \mathbf{Q}/\mathbf{Z}
\end{displaymath}
and defined by: 
\begin{equation}
\label{eq:defquad}
\phi_{f,w}(\bar{x})=\frac{1}{2}\left(x_{\mathbf{Q}}(\tilde{x})-x(w)\right)
\quad mod \quad 1.
\end{equation}
$\phi_{f,w}$ is determined by the modulo $2 F$ class of $w$. The triple
$(F,f,w)$ is said to be a \emph{presentation} of the quadratic
 form $\left(K_{f},\phi_{f,w}\right)$.\\
Any linking pairing and any nondegenerate quadratic form
admit such presentations with $f$ nondegenerate 
(see \cite[Theorem (6)]{WalI}).

Given an arbitrary quadratic form on a finite Abelian group $(G,q)$,
we can calculate its \emph{Gauss sum}:
\begin{displaymath}
\gamma(G,q)=\frac{1}{\sqrt{|ker(b_{q})| \cdot |G|}}\cdot
            \sum_{\substack{x\in G}} e^{2i \pi q(x)} \in \mathbf{C}.
\end{displaymath}
This complex number is a $8^{\textrm{e}}$-root of unity or
is $0$ (only if $q$ is degenerate).\\ 
We then define the corresponding \emph{Gauss-Brown invariant}
$B(G,q)\in \overline{\mathbf{Z}_{8}}=\mathbf{Z}_{8} \cup \{ \infty \}$ 
by the formula:
\begin{equation}
\gamma(G,q)=e^{\frac{2i\pi}{8}\cdot B(G,q)} \in \mathbf{C},
\end{equation} 
using the convention $e^{2i\pi\cdot \infty}=0$.\\
If $(G,q)$ admits a triple $(F,f,w)$ as a presentation,
a useful formula of Van der Blij states that\footnote{A sketch 
of proof of this formula can be found in \cite{VDB},
in case when $f$ is nondegenerate.\\ 
Milnor and Husemoller have included a detailed proof 
in \cite[Appendix 4]{MH}, for $f$ nondegenerate
and $w=0$. The general case can be reduced to this special case.}: 
\begin{equation}
\label{eq:VDB}
B(K_{f},\phi_{f,w})=sgn(f)-f(w,w) \quad mod  \quad 8.
\end{equation} 

For $(G,b)$ a linking pairing, denote by $Quad(G,b)$ the set of quadratic
forms with $b$ as associated linking pairing, and denote by
$T_{2}(G)$ the subgroup of elements of $G$ of order at most $2$. Using
the nondegenerativity of $b$, we easily obtain:
\begin{lemma}
\label{lem:Quad_is_affine}
The set $Quad(G,b)$ is an affine space over $T_{2}(G)$, 
with action defined by:
\begin{displaymath}
\forall x \in T_{2}(G), \ \forall q \in Quad(G,b),  \
x\cdot q:=q+b(x,-) .
\end{displaymath}             
\end{lemma}
The following lemma says how the Gauss-Brown 
invariant behaves under this action.
\begin{lemma}
\label{lem:how_Gauss-Brown_behaves}
Let $q\in Quad(G,b)$ and let $x\in T_{2}(G)$. Then,
\begin{displaymath}
\gamma(G,x\cdot q)=e^{-2i\pi q(x)}\cdot \gamma(G,q).
\end{displaymath}
\end{lemma}
\begin{proof}
Since $(x\cdot q)(y)=q(y)+b(x,y)=q(x+y)-q(x)$, we have:
\begin{displaymath}
\begin{array}{rcl}
\sqrt{|G|}\cdot \gamma(G,x\cdot q) & = & \sum_{y\in G} e^{2i \pi (x\cdot q)(y)} \\
                   & = & e^{-2i\pi q(x)}\cdot \sum_{y\in G} e^{2i \pi q(x+y)}\\
& = & e^{-2i\pi q(x)} \cdot \sqrt{|G|} \cdot \gamma(G,q)
\end{array}
\end{displaymath}
\end{proof} 
\subsection{Isomorphism classes of nondegenerate quadratic forms}
\label{subsec:betq}
We now want to prove the following result:
\begin{theorem}
\label{th:betq}
Let $(G,q)$ and $(G',q')$ be two nondegenerate quadratic forms on finite 
Abelian groups. We denote by $b$ and $b'$ the linking pairings 
going respectively with them.
The following two assertions are equivalent:
\begin{enumerate}
\item $(G,q)\simeq (G',q')$,
\item $(G,b) \simeq (G',b')$ and $B(q)=B(q') \in \mathbf{Z}_{8}$,
\end{enumerate}
where $B(-)$ denotes the Gauss-Brown invariant of quadratic forms.
\end{theorem}
Let $(G,q)$ (respectively $(G',q')$) be a nondegenerate quadratic 
form over the linking pairing $(G,b)$ (respectively $(G',b')$).
Recall that the relation between $b$ and $q$ is:
\begin{equation}
\label{eq:betq}
b(x,y)=q(x+y)-q(x)-q(y).
\end{equation}
From this formula and the definition of $B(-)$,
``$1\Rightarrow 2$'' of Theorem \ref{th:betq} is obvious.

Suppose momentarily that $G$ is a $p$-group, with $p$ an odd prime.
Then, equation $2q(x)=b(x,x)$ makes $b$ determine $q$,
for if $q''$ is another quadratic form over $b$, then $(q-q'')$ is
an order at most $2$ element of 
$Hom\left(G,\mathbf{Q}/\mathbf{Z}\right)\simeq G$, 
and so vanishes.
So, if $G$ and $G'$ are both $p$-groups, then 
$b\simeq b'$ implies $q\simeq q'$.

Come back now to the general case and
suppose that condition 2 is satisfied. 
$(G,b)$ splits along its $p$-primary components $G_{p}$:
\begin{displaymath}
(G,b)= \bigoplus_{p,\textrm{ prime }} 
\left(G_{p},b|_{G_{p}\times G_{p}}\right),
\end{displaymath} 
and, according to formula (\ref{eq:betq}), the same holds for $q$.
The given isomorphism between $(G,b)$ and $(G',b')$, 
induces then for each prime $p$
an isomorphism between $\left(G_{p},b|_{G_{p}\times G_{p}}\right)$ and 
$\left(G'_{p},b'|_{G'_{p}\times G'_{p}}\right)$. 
From the above lines, we deduce that,
for $p$ odd, $\left(G_{p},q|_{G_{p}}\right)$ 
and $\left(G'_{p},q'|_{G'_{p}}\right)$ are isomorphic.\\
In particular, $B\left(q|_{G_{p}}\right)=
B\left(q'|_{G'_{p}}\right)$ for $p$ odd, and so,
by additivity of the Gauss-Brown invariant, this is also true for $p=2$.
Consequently, it is enough to prove Theorem \ref{th:betq}
when $G$ and $G'$ are $2$-groups.

In the sequel, we recall a construction due to Wall 
(see \cite[\S 6]{WalI}),
establishing a one to one correspondence (up to isomorphism) between 
nondegenerate quadratic forms on $2$-groups and linking pairings 
on $2$-groups without direct summand
of order two. Next, we give a brief review of Kawauchi and 
Kojima classification of linking pairings on $2$-groups. We will finally
end the proof of Theorem \ref{th:betq}.
\subsubsection*{Wall construction}
A linking pairing $(G',b')$ will be said here to be 
\emph{special} if $G'$ is a finite $2$-group without
direct summand of order two.\\
Let give us a special linking pairing $(G',b')$. 
Set $G:=G'/ T_{2}(G')$ where $T_{2}(G')$ is the subgroup of elements of 
order at most $2$, and denote by $\pi$ the canonical projection $G'\rTo G$.
Note that $G$ can be any $2$-group.
Define now $q: G \rTo \mathbf{Q}/\mathbf{Z}$ by: 
\begin{displaymath} 
\forall x=\pi(x')\in G, \quad q(x)=b'(x',x').
\end{displaymath}
The quantity $q(x)$ is well-defined because of the special feature of the
group $G'$.\\
Then, $q$ is easily seen to be quadratic and nondegenerate, 
with associated linking pairing $b:G\times G 
\rTo \mathbf{Q}/\mathbf{Z}$ defined by:
\begin{equation}
\label{eq:betbprime}
\forall x=\pi(x'),\ \forall y=\pi(y')\in G, \ b(x,y)=2\cdot b'(x',y').
\end{equation}
Let us denote by $\Psi$ the construction $(G',b')\rMapsto (G,q)$. Wall
showed this to be surjective onto the set
of nondegenerate quadratic forms on $2$-groups. He also proved that 
if $(G',b'_{1})$ and $(G',b'_{2})$
give rise to the same $(G,q)$ by $\Psi$, then they have to be isomorphic 
(see \cite[Theorem 5]{WalI}).

As a consequence, the classification, up to isomorphism,
of nondegenerate quadratic forms on $2$-groups 
is reduced to that of special linking pairings.
\subsubsection*{Kawauchi-Kojima classification of linking pairings on 2-groups}
Let $(G,b)$ be a linking pairing on a finite $2$-group. Find a cyclic
decomposition of $G$:
\begin{displaymath}
G\simeq \bigoplus_{k\geq 1} \left(\mathbf{Z}_{2^{k}}\right)^{r_{k}}.
\end{displaymath}
The natural numbers $r_{k}=r_{k}(b)$ are group invariants of $G$.
The very next construction is due to Wall (see \cite[\S 5]{WalI}).

Denote by $\overline{G}_{k}$ the subgroup of $G$ 
of elements of order at most $2^{i}$ 
for $i\leq k$, and set:
\begin{displaymath}
\tilde{G}_{k}=\frac{\overline{G}_{k}}
                   {\overline{G}_{k-1}+2\cdot \overline{G}_{k+1}}.
\end{displaymath}
The group $\tilde{G}_{k}$ is clearly a $\mathbf{Z}_{2}$-vector space of 
rank $r_{k}$.
Let also $\tilde{b}_{k}: \tilde{G}_{k} \times \tilde{G}_{k}
\rTo \mathbf{Q}/\mathbf{Z}$ be defined by:
\begin{displaymath}
\forall x,y \in \overline{G}_{k},\quad
\tilde{b}_{k}(\tilde{x},\tilde{y})=2^{k-1}\cdot b(x,y).
\end{displaymath}
The form $\tilde{b}_{k}$ was shown by Wall to be nondegenerate
(see \cite[Lemma 8]{WalI}).

Consider now the map $\tilde{G}_{k}\rTo \mathbf{Q}/\mathbf{Z}$
sending $\tilde{x}$ to $\tilde{b}_{k}(\tilde{x},\tilde{x})$. It is additive, 
and so we can define an element $\tilde{c}_{k}=\tilde{c}_{k}(b)$ in 
$\tilde{G}_{k}$ by the equation:
\begin{equation}
\forall \tilde{x} \in \tilde{G}_{k},\quad
\tilde{b}_{k}(\tilde{x},\tilde{c}_{k})=\tilde{b}_{k}(\tilde{x},\tilde{x}).
\end{equation}
When $\tilde{c}_{k}=0$, the following map can be defined:
\begin{displaymath}
\left\{\begin{array}{rrcl}
\frac{G}{\overline{G}_{k}} & \rTo^{q_k} & 
\frac{\mathbf{Q}}{\mathbf{Z}} \\
\bar{x} & \rMapsto & 2^{k-1}\cdot b(x,x).
\end{array}\right.
\end{displaymath}
The form $q_{k}=q_{k}(b)$ can be verified to be quadratic and nondegenerate. 
In particular, its Gauss-Brown invariant $B(q_{k})$ is not equal to $\infty$.\\
Kawauchi and Kojima defined $\sigma_{k}
=\sigma_{k}(b) \in \overline{\mathbf{Z}_{8}}$ by:
\begin{displaymath}
\sigma_{k}=\left\{ \begin{array}{ll}
                       B(q_{k}) & \textrm{if } \tilde{c}_{k}=0 ,\\
                       \infty & \textrm{otherwise,}
                       \end{array}
              \right.
\end{displaymath}              
and showed the following theorem (see \cite[Theorem 4.1]{KK}).
\begin{theorem}[Kawauchi-Kojima]
If $(G,b)$ is a linking pairing on a finite $2$-group,
its isomorphism class is determined by the invariant family
$\left(r_{k}(b),\sigma_{k}(b)\right)_{k\geq 1}$.
\end{theorem}

\begin{proof}[End of proof of Theorem \ref{th:betq}]
Let $(G,q)$ be a quadratic form on a finite $2$-group $G$, going with
a linking pairing $b$. 
Let also $(G',b')$ be a special linking pairing, 
giving rise to $(G,q)$ by Wall construction $\Psi$.

We want to compare the invariants $r_{k}$ and $\sigma_{k}$ of $b$ and $b'$,
in order to quantify how much $(G,b)$ determines $(G',b')$, and so 
$(G,q)$, up to isomorphism.\\
Recall that $G=G'/ \overline{G'}_{1}$. Denote by $\pi:G'
\rTo G$ the canonical projection. The map $\pi$
induces a morphism from $\overline{G'}_{k+1}$ onto $\overline{G}_{k}$ 
with kernel $\overline{G'}_{1}$. So, $\pi$ induces a natural isomorphism:
\begin{diagram}
\tilde{G'}_{k+1} & \rTo^{\tilde{\pi}_{k}}_{\simeq} & \tilde{G}_{k}.
\end{diagram}
In particular, the $\mathbf{Z}_{2}$-vector spaces $\tilde{G'}_{k+1}$ and
$\tilde{G}_{k}$ have the same rank. So,
\begin{equation}
\label{eq:r}
\forall k\geq 1, \quad r_{k+1}(G')=r_{k}(G).
\end{equation}
Besides, since $G'$ is special, we have: 
\begin{equation}
\label{eq:r1}
r_{1}(G')=0.
\end{equation}
The isomorphism  $\tilde{\pi}_{k}$ makes $\tilde{b'}_{k+1}$ and
$\tilde{b}_{k}$ commute because of equation (\ref{eq:betbprime}).
As a consequence, $\tilde{\pi}_{k}$ sends $\tilde{c}_{k+1}(b')$ to 
$\tilde{c}_{k}(b)$. 
Furthermore, when these (simultaneously) vanish, the natural isomorphism 
between $G'/\overline{G'}_{k+1}$ and $G/\overline{G}_{k}$ induced by $\pi$, 
make $q_{k+1}(b')$ and $q_{k}(b)$ commute (because
of equation (\ref{eq:betbprime})). As a consequence, these two
quadratic forms will have the same Gauss-Brown invariant. So, to sum up,
\begin{equation}
\label{eq:sigma}
\forall k\geq 1, \quad
\sigma_{k+1}(b')=\sigma_{k}(b)\in \overline{\mathbf{Z}_{8}}.
\end{equation}
Since $\tilde{G'}_{1}=0$, $\tilde{c}_{1}(b')$ vanishes.
It remains to be noticed that $q_{1}(b')$ is nothing but $q$. Thus,
\begin{equation}
\label{eq:sigma1}
\sigma_{1}(b')=B(q).
\end{equation}

Now, from equations (\ref{eq:r}), (\ref{eq:r1}),
(\ref{eq:sigma}), (\ref{eq:sigma1}) and Kawauchi-Kojima theorem,
we see that $(G,b)$ together with $B(q)$ determine $(G,q)$ up to isomorphism.
What has been remaining to be proved for Theorem \ref{th:betq}, 
then follows.
\end{proof}

We now give a result of Durfee (see \cite[Corollary 3.9]{Dur}) as a corollary
of Theorem~\ref{th:betq}.
\begin{corollary}
\label{cor:noloworder}
Let $b:G\times G\rTo \mathbf{Q}/\mathbf{Z}$ be a linking pairing
on a finite Abelian group $G$ \emph{without cyclic direct summand
of order $2$ or $4$}.
Then, $\forall q,q' \in Quad(G,b),\ q\simeq q'$.
\end{corollary}
\begin{proof}  Take some quadratic forms $q$ and $q'$ over $b$, and
let $x\in T_{2}(G)$ be such that $q'=x\cdot q$ (see  
Lemma \ref{lem:Quad_is_affine}). By the hypothesis on $G$,
there exits some $x_{0} \in G$ such that $x=4x_{0}$, and so,
by homogeneity of $q$, we have:
\begin{displaymath}
q(x)=q(4x_{0})=4^{2}q(x_{0}).
\end{displaymath}
But since $x_{0}$ is then of order at most $8$, $q(x_{0})$ has to be of order
at most $2\cdot 8=16$ (see \cite[Lemma 1.12]{DelII}). It follows
that $q(x)=0$ and so, by Lemma \ref{lem:how_Gauss-Brown_behaves},
we obtain that $B(G,q)=B(G,q')$. Theorem \ref{th:betq} allows us
to conclude.
\end{proof}
\section{The quadratic form $\phi_{M,\sigma}$}
\label{sec:topquad}
In this section, when not specified, integer coefficients are assumed.
\subsection{Turaev $4$-dimensional definition of $\phi_{M,\sigma}$}
\label{subsec:dim4}
Let $M$ be a  connected closed oriented $3$-manifold, and let 
$\psi: V_{L} \rTo M$ be a surgery presentation on $\mathbf{S}^{3}$
given by an ordered oriented framed link $L$ 
(see the beginning of \S\ref{subsec:surgery_presentations}).\\
We  use notations and apply constructions 
of \S\ref{subsec:presentations} to 
$F=H_{2}(W_{L})$, taking for $f$ the intersection form of $W_{L}$.
Recall that the matrix of $f$ relative 
to the preferred basis of $F$
is $B_{L}$, the linking matrix of $L$.\\
The composite: 
\begin{diagram}
H_{2}(W_{L}) & \rTo^{ad_{f}} & Hom(H_{2}(W_{L}),\mathbf{Z}) &
\lTo_{\simeq} &  H^{2}(W_{L}) &\rTo^{P}_{\simeq} & 
H_{2}(W_{L},\partial W_{L})
\end{diagram}
is equal to 
$i_{*}:H_{2}(W_{L}) \rTo 
H_{2}(W_{L},\partial W_{L})$,
induced by inclusion.
Since $Coker(i_{*})=H_{1}(V_{L})$,
we obtain the following isomorphism $r$:
\begin{diagram}
K_{f} & \rTo_\simeq &  T(Coker( i_{*})) & \rTo_\simeq^{\psi_{*}} &
TH_{1}(M)\\
\dDash & & & & \uDashto \\
\HmeetV & & \rDash~{r} & & \HmeetV
\end{diagram}
In fact, it is well-known that the above $(F,-f)$ is 
\emph{via} $r$ a presentation
of the torsion linking form $\big( TH_{1}(M), \lambda_{M} \big)$ of $M$,
the definition of which we now recall:
\begin{quote}
Let $x,x' \in TH_{1}(M)$ be respectively realized by oriented 
disjoint knots $K,K'$ in $M$.
Let $c' \in \mathbf{N}$ be such that $c'\cdot x'=0$.
Pick a $c'$-times connected sum of $K'$. 
We obtain a null-homologous knot in $M$ for which we
can thus find a Seifert surface $S'$ in general position with $K$. Then:
\begin{displaymath}
\lambda_{M}(x,x') =\frac{1}{c'} K \bullet S'\in \mathbf{Q}/\mathbf{Z} ,
\end{displaymath}
where $\bullet$ is the intersection form of $M$.
\end{quote}
Now to each $\sigma \in Spin(M)$ is associated a characteristic
solution of $-B_{L}$ or, alternatively, a Wu class (modulo $2 F$) of $-f$,
denoted by $w_{\sigma}$. Then, Turaev defined:
\begin{definition}
\label{def:dim4}
The \emph{quadratic form of the spin $3$-manifold $(M,\sigma)$}:
\begin{diagram}
TH_{1}(M) & \rTo^{\phi_{M,\sigma}} & \mathbf{Q} / \mathbf{Z}
\end{diagram}
is defined to be $\phi_{-f,w_{\sigma}} \circ r^{-1}$.  
\end{definition} 
We still have to verify that $\phi_{M,\sigma}$ does not depend on the choice
of the surgery presentation.\\
Let $\psi':V_{L'} \rTo M$ be another one.
Let $s \in (\mathbf{Z}_{2})^{l}$ (resp. $s'$) 
be the characteristic solution of 
$B_{L}$ (resp. $B_{L'}$) corresponding to the spin structure
 $\psi^{*}(\sigma)$ of $V_{L}$ (resp. $(\psi')^{*}(\sigma)$ of $V_{L'}$).\\
According to the refined Kirby theorem (see \cite{Bl}), 
there  exists a sequence of spin
Kirby moves from $(L,s)$ to $(L',s')$, inducing
a spin-diffeomorphism from $\left(V_{L},\psi^{*}(\sigma)\right)$ to 
$\left(V_{L'},(\psi')^{*}(\sigma)\right)$ 
isotopic to $(\psi')^{-1} \circ \psi$.
These Kirby moves induce a path 
$(F,f,w_{\sigma}) \leadsto (F',f',w'_{\sigma})$ whose elementary steps are:
\begin{displaymath}
\left\{\begin{array}{ll}
(F,f,w) \rMapsto \left(F,{}^{t}S f S,S^{-1}(w)\right) & 
\textrm{with } S\in Aut(F),\\
(F,f,w) \rMapsto \left(F\oplus \mathbf{Z}, 
f\oplus (\pm 1),w\oplus(1)\right). &
\end{array} \right.
\end{displaymath}
As a consequence, this path induces an isomorphism $t$ from
$\left(K_{f},\phi_{f,w_{\sigma}}\right)$ to 
$\left(K_{f'},\phi_{f',w'_{\sigma}}\right)$ making
the following diagram commutative:
\begin{diagram}
K_{f} & \rTo^{t} & K_{f'} \\
& \rdTo<r & \dTo>{r'}\\
& & TH_{1}(M)
\end{diagram}
The well-definition of $\phi_{M,\sigma}$ then follows. 
\subsection{An intrinsic definition for $\phi_{M,\sigma}$}
\label{subsec:dim3}
Let $(M,\sigma)$ be a closed spin $3$-manifold and
let $K$ be a smooth oriented knot in $M$.\\
Each parallel $l$ of $K$ defines a trivialization of the normal
bundle of $K$ in $M$, and so allows us to restrict $\sigma$
to a spin structure on $K\cong \mathbf{S}^{1}$. We define:
\begin{displaymath}
\sigma(K,l)=\left\{ \begin{array}{ll}
                    0 & \textrm{ if the induced spin structure on }
                        \mathbf{S}^{1}
                        \textrm{ is the ``bounding'' one,}\\
                    1 & \textrm{ otherwise.}
                    \end{array} \right.
\end{displaymath}
$\sigma(K,l)$ is a $\mathbf{Z}_{2}$-valued invariant of framed
knots in $M$.
\begin{lemma}
\label{lem:framed}
Let $(M,\sigma)$ be a closed spin $3$-manifold. Then,
for each oriented smooth knot $K$ in $M$ 
with $l$ as a parallel and meridian $\mu$,
\begin{displaymath}
\sigma(K,l+\mu)=\sigma(K,l)+1 \in \mathbf{Z}_{2}.
\end{displaymath}
\end{lemma}
\begin{proof}  Let $S$ denote the boundary of a regular neighbourhood
of $K$ in $M$. The normal bundle of $S$ in $M$ is naturally trivialized,
so $S$ inherits from $(M,\sigma)$ a spin structure. Let $q$ be
the quadratic form associated 
to the spin smooth surface $\left(S,\sigma|_{S}\right)$
as defined by Johnson in \cite{Joh}. The following identity then holds
for each parallel $l$:
\begin{displaymath}
\sigma(K,l)=q\left([l]\right) ,
\end{displaymath}
when $l$ is thought of as a curve on $S$. Since $q$ is quadratic
with respect to the modulo $2$ intersection form $\bullet$ on $S$, we have:
\begin{displaymath}
\begin{array}{rcl}
\sigma(K,l+\mu)&=&q\left([l]+[\mu]\right) \\
               &=&q\left([l]\right)+q\left([\mu]\right)+[l]\bullet [\mu] \\
               &=& q\left([l]\right) +1 \\
               &=& \sigma(K,l) +1.
\end{array}
\end{displaymath}
\end{proof}
We now recall the definition of the \emph{framing number}
$Fr(K,l) \in \mathbf{Q}$ of a  rationally nulhomologous oriented framed knot
$(K,l)$ in a closed oriented $3$-manifold $M$:
\begin{quote}                          
Choose $c\in \mathbf{N}^{*}$ such that $c\cdot[K]=0\in H_{1}(M)$.
Pick a $c$-times connected sum of $K$. We obtain a null-homologous knot in
$M$ for which we can thus pick a Seifert surface $S$ in
general position with the knot $l$. Then:
\begin{displaymath}
Fr(K,l)=\frac{1}{c}\cdot l\bullet S \in \mathbf{Q} .
\end{displaymath}                                 
\end{quote}
\begin{lemma}
\label{lem:dim3}
Let $(M,\sigma)$ be a closed spin $3$-manifold and $x\in TH_{1}(M)$.
Choose a smooth oriented knot $K$ in $M$ representative for $x$, and pick
a parallel $l$ for $K$ satisfying $\sigma(K,l)=0\in \mathbf{Z}_{2}$. Then,
\begin{equation}
\label{eq:dim3}
\phi_{M,\sigma}(x)=\frac{1}{2}\cdot Fr(K,l) 
\in \frac{\mathbf{Q}}{\mathbf{Z}}.
\end{equation}
\end{lemma}
Note that, according to Lemma \ref{lem:framed},
the rhs of (\ref{eq:dim3}) is an invariant of the oriented knot $K$ (it does
not depend on the choice of $l$ satisfying the above condition).
This lemma claims that it only depends on the homology class $x$
of $K$, and \emph{gives a $3$-dimensional 
definition for the quadratic form $\phi_{M,\sigma}$}.
\begin{remark}
From this lemma, we can see that $\phi_{M,\sigma}$ coincides with the
quadratic form defined by Lannes and Latour in \cite{LL} when specialized
to our case (see also \cite{MS}).
\end{remark}
\begin{proof}[Proof of Lemma \ref{lem:dim3}] 
Consider the $4$-manifold $W_{1}$ obtained from $M\times[0,1]$
by attaching a $2$-handle to $M\times 1$ along $(K,l)$. Identify
$M$ with $M\times 0$. Since $\sigma(K,l)=0$, $\sigma$ extends in
a unique way to a spin structure $\sigma_{1}$ of $W_{1}$.\\
$(W_{1},\sigma_{1})$ is then a spin cobordism between $(M,\sigma)$ and
$(-M',-\sigma')$, where $M'$ is the closed oriented
$3$-manifold obtained from $M$ by the corresponding surgery, 
and where $-\sigma'$ is the restriction of $\sigma_{1}$ to $-M'$.\\
Note also that the core of the $2$-handle is a $2$-disc $D$ in $W_{1}$
with boundary $K$ in $M$, and 
whose normal bundle can be trivialized in accordance
with the trivialization of the normal bundle of $K$ in $M$ given by $l$.
The framed knot $(K,l)$ will briefly be said to \emph{have property 
$(\mathcal{D})$} in $W_{1}$.\\
According to Kaplan Theorem (see \cite{Kap}), the spin $3$-manifold
$(M',\sigma')$ admits an \emph{even} surgery presentation in $\mathbf{S}^{3}$
(\emph{i.e.} the linking matrix is even and its characteristic solution
corresponding to $\sigma'$ is the trivial one). Denote by $W_{2}$ the trace of 
the surgery and by $\sigma_{2}$ the unique extension of $\sigma'$
to the whole of $W_{2}$.\\
By glueing $(W_{1},\sigma_{1})$ to  $(W_{2},\sigma_{2})$ along
$(M',\sigma')$, we obtain a spin $4$-manifold $(W,\sigma)$ with
boundary $(M,\sigma)$. The $2$-handle from $M$ to $M'$ can be reversed.
After a rearrangement, the $4$-manifold $W$ appears as 
$\mathbf{B}^{4}$ to which
have been simultaneously attached some $2$-handles (one more than $W_{2}$),
with boundary $M$, and to which $\sigma$ can be extended. So,
$W$ is the trace of an even surgery presentation.

\emph{As a summary}, we have found so for an even surgery presentation $(L,0)$
for $(M,\sigma)$ such that $(K,l)$ has property ($\mathcal{D}$) in
$W_{L}$.

Let us work with this surgery presentation of $(M,\sigma)$.
Notations of \S\ref{subsec:dim4} will be used: $F=H_{2}(W_{L})$,
$f$ stands for the intersection form of $W_{L}$ and so on.\\
Let the $2$-disc $D$ give an element of $H_{2}(W_{L},\partial W_{L})$.
The latter is identified with an element $d$ of $Hom(F,\mathbf{Z})$.
Recall from the definition of $\phi_{M,\sigma}$ 
that, in this even case,
\begin{equation}
\label{eq:tilde}
\phi_{M,\sigma}(x)=-\frac{1}{2}d_{\mathbf{Q}}(\tilde{d}) ,
\end{equation}
where $d_{\mathbf{Q}}$ is the rational extension of $d$ and where 
$\tilde{d}\in F\otimes \mathbf{Q}$ is such that $a_{f}(\tilde{d})=d$.\\
Let $c \in \mathbf{N}$ such that $c\cdot x=0\in H_{1}(M)$. Then, there exists
$y\in F$ such that $ad_{f}(y)=c\cdot d$. So 
$\frac{1}{c}\cdot y \in F\otimes \mathbf{Q}$
works as a $\tilde{d}$. Equation (\ref{eq:tilde}) can be rewritten as:
\begin{equation}
\label{eq:notilde}
\phi_{M,\sigma}(x)=-\frac{1}{2c^{2}}f(y,y) .
\end{equation}
When $y$ is seen as belonging to $H_{2}(W_{L})$, the integer $f(y,y)$ is
equal to $Y\bullet Y'$, where $Y$ and $Y'$ are $2$-cycles representatives
for $y$ in transverse position in $W_{L}$.
By means of a ``collar'' trick appearing in \cite{DelII},
we will be able to give
examples of such $Y$ and $Y'$.\\
We add a collar $M\times [0,1]$ to $W_{L}$ 
such that $M\times 0$ is identified with $M$.
Let $S$ be a Seifert surface for $c\cdot K$ 
in $M$ in transverse position with $l$,
and $S'$ be a Seifert surface for $(c\cdot l)\times 1$ in $M\times 1$.
Because of the property $(\mathcal{D})$, $D$ can be pushed off to a disc $D'$
in such a way that $\partial D'=l$ and $D\cap D'=\varnothing$.
Figure \ref{fig:collar} is a good summary.
\begin{figure}[!h]
\begin{center}
\includegraphics[width=6cm,height=5cm]{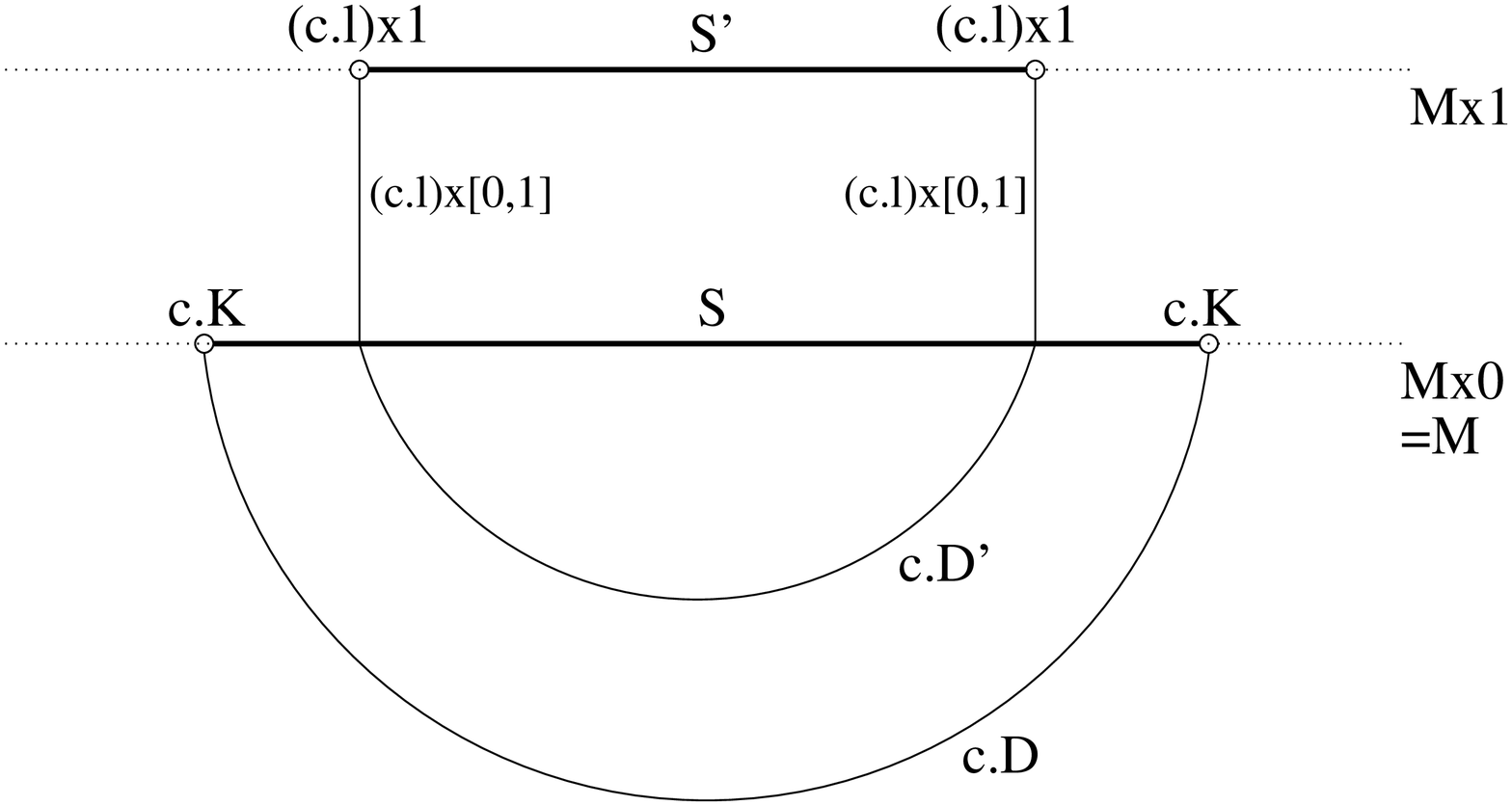}
\caption{Calculating $Y\bullet Y'$.}
\label{fig:collar}
\end{center}
\end{figure}
We define $Y=c.D-S$ and $Y'=
\left(c.D'+(c.l)\times [0,1]\right)-S'$. Then:
\begin{equation}
\label{eq:WtoM}
f(y,y) = Y\bullet Y'
       = (-S)\bullet(c.D')
       = -c\cdot S\bullet l .
\end{equation}
The last $\bullet$ in (\ref{eq:WtoM}) is intersection in $M$.
The lemma follows from (\ref{eq:WtoM}), (\ref{eq:notilde})
and the  definition of a framing number.
\end{proof} 
\subsection{Properties of  $\phi_{M,\sigma}$}
\label{subsec:properties}
The algebraic results of \S\ref{sec:algquad} have a topological
meaning. First, it has been shown by Turaev in \cite[Theorem V]{Tur}:
\begin{lemma}[Turaev]
\label{lem:BetR}
For each $\sigma \in Spin(M)$, 
$B(\phi_{M,\sigma})=-R_{M}(\sigma) \quad mod \quad 8$,\\
where $R_{M}$ is the Rochlin function of the $3$-manifold $M$.
\end{lemma}
\begin{proof}  Find an \emph{even} surgery presentation
of the spin $3$-manifold $(M,\sigma)$. So
we are led to apply formula (\ref{eq:VDB}) with $w=0$. 
\end{proof}
So, in view of Lemma \ref{lem:BetR}, the topological
translations of Theorem \ref{th:betq} 
and its Corollary \ref{cor:noloworder} are respectively:
\begin{proposition}
\label{prop:lambdaetphi}
Let $(M,\sigma)$ and $(M',\sigma')$ be connected closed  spin $3$-manifolds.\\
The following two assertions are equivalent:
\begin{enumerate}
\item their quadratic forms $\phi_{M,\sigma}$ and $\phi_{M',\sigma'}$
      are isomorphic,
\item their linking forms $\lambda_{M}$ and $\lambda_{M'}$ are isomorphic,
       and $R_{M}(\sigma)=R_{M'}(\sigma')$ modulo $8$.
\end{enumerate}
\end{proposition}
\begin{corollary}
\label{cor:top_noloworder}
Let $M$ be a connected  closed oriented $3$-manifold, such that $H_{1}(M)$
does not admit $\mathbf{Z}_{2}$ nor $\mathbf{Z}_{4}$ as a direct summand.
Then, all of its quadratic forms are isomorphic one to another.
\end{corollary}
\section{Proof of refined Matveev theorem}
\label{sec:proof}
Part of the work has already been done in previous sections. 
First, ``$1\Longrightarrow 3$'' 
follows from Corollary \ref{cor:Rochlin} 
and from the easy part of (unspun) Matveev theorem: a $Y$-surgery preserves
homology and torsion linking forms. This can be verified seeing $Y$-surgery
as a $\mathcal{V}_{h}$-surgery (where $h$ is the Borromean diffeomorphism of
Remark \ref{rem:fundamental}), using a Mayer-Vietoris argument and 
the fact that $h$ belongs to the Torelli group.
Second, ``$3\Longrightarrow 2$'' follows from Proposition                                              
\ref{prop:lambdaetphi}. What remains to be proved is then
``$2\Longrightarrow 1$''.\\
                                     
We start by recalling an algebraic result of Durfee
(see \cite{Dur}) about
even symmetric bilinear forms on finitely generated free Abelian groups.\\ 
Let $(\mathbf{Z}^{2},h)$ and 
$(\mathbf{Z}^{8},\gamma_{8})$ be the unimodular even forms
whose matrices are respectively:
\begin{displaymath}
H=\left( \begin{array}{cc}
          0 & 1 \\
          1 & 0
          \end{array} \right)
\qquad \textrm{and} \qquad
\Gamma_{8}=\left( \begin{array}{cccccccc}
                   2 & 1 &&&&&&\\
                   1 & 2 & 1 &&&&&\\
                   & 1 & 2 & 1 &&&&\\
                   & & 1 & 2& 1&&&\\
                   &&& 1 & 2 & 1& 0& 1 \\
                   &&&&1&2&1&0\\
                   &&&&0&1&2&0\\
                   &&&&1&0&0&2
                   \end{array} \right).
\end{displaymath}
\begin{definition}
\label{def:stable}
Let $(F_{1},f_{1})$ and $(F_{2},f_{2})$ be two symmetric even bilinear forms
on finitely generated free Abelian groups.
They are said to be \emph{stably equivalent}
if  they become isomorphic after some stabilizations 
with unimodular  symmetric even bilinear forms.
\end{definition}                    
Note that a unimodular even form becomes indefinite after taking direct
sum with $h$. Recall that every even unimodular indefinite form splits as 
a direct sum of $h$ and $\gamma_{8}$ (see for example \cite[p.26]{Kir}). 
Thus, in Definition \ref{def:stable},
$h$ and $\gamma_{8}$ suffice as unimodular
even forms to stabilize.
 
For $(F,f)$ an even symmetric bilinear form
on a finitely generated free Abelian group, we will shortly denote
by $\phi_{f}$ the quadratic form $\phi_{f,0}$ corresponding to the zero
Wu class.
\begin{proposition}[Durfee]
\label{prop:Durfee}
Let $(F_{1},f_{1})$ and $(F_{2},f_{2})$ be two symmetric even bilinear forms
on finitely generated free Abelian groups.
Then, the following two assertions are equivalent:
\begin{enumerate}
\item $(F_{1},f_{1})$ and $(F_{2},f_{2})$ are stably equivalent,
\item $Ker\left(ad_{f_{1}}\right) \simeq 
Ker\left(ad_{f_{2}}\right)$ and 
$\phi_{f_{1}} \simeq \phi_{f_{2}}$.
\end{enumerate}
\end{proposition}
\begin{proof}  Implication ``1$\Rightarrow$2'' 
is obvious since $\gamma_{8}$
and $h$ are both unimodular.\\
Now suppose that condition 2 is satisfied. For each $i\in\{1,2\}$,
there exists a nondegenerate symmetric bilinear form 
$\left(\tilde{F}_{i},\tilde{f}_{i}\right)$ such that:
\begin{displaymath}
(F_{i},f_{i})\simeq \left(\tilde{F}_{i},\tilde{f}_{i}\right)
\oplus (\mathbf{Z}^{r},O_{r})
\end{displaymath}
where $\tilde{F}_{i}=F_{i}/Ker\left(ad_{f_{i}}\right)$ 
(note that it is free),
$r=rk\left(Ker\left(ad_{f_{1}}\right)\right)
=rk\left(Ker\left(ad_{f_{2}}\right)\right)$,
and $O_{r}$ is the zero form.
The form $\tilde{f}_{i}$ is still even and 
$\phi_{\tilde{f}_{1}}$ and $\phi_{\tilde{f}_{2}}$ are still isomorphic.\\
Consequently, without loss of generality, we can assume both $f_{i}$
to be nondegenerate. But this case was treated by Durfee in
\cite[Corollary 4.2.(ii)]{Dur}\footnote{An easier 
and direct proof of this result was given 
by Wall in \cite[Corollary 1]{WalII}.}.
\end{proof}
Let $(M,\sigma)$ and $(M',\sigma')$ be connected closed  
spin $3$-manifolds such that
$\beta_{1}(M)=\beta_{1}(M')$ and $\phi_{M,\sigma} \simeq \phi_{M',\sigma'}$.\\
Suppose given for them
some even surgery presentations $(L,0)$ and $(L',0)$ with 
respective linking matrices $B$ and $B'$.
According to Proposition \ref{prop:Durfee}, there exists a unimodular 
integer matrix $P$ satisfying for some stabilizations:
\begin{displaymath}
{}^{t}P\cdot (B\oplus H \cdots H \oplus 
\Gamma_{8} \cdots \Gamma_{8})\cdot P
=B'\oplus H \cdots H \oplus \Gamma_{8} \cdots \Gamma_{8} .
\end{displaymath}
We have the following geometric realizations of algebraic operations:
\begin{enumerate}
\item stabilizations by $H$ correspond to connected sums with 
$\mathbf{S}^{3}$ surgery presented on the zero-framed Hopf link,
\item a stabilization by $\Gamma_{8}$ is concrete when thought of as a 
connected sum with the Poincar\'e sphere surgery presented on an 
appropriate height-component link as in \cite[Figure 5.3, p.15]{Kir},
\item congruence by $P$ can be realized by some spin Kirby moves
(handle-slidings and changes of orientation of components of $L$).
\end{enumerate}
The Poincar\'e sphere can also be obtained by surgery along a $(+1)$-framed
trefoil knot (\cite[Figure 5.3, p.15]{Kir}), which can be obtained from the
$(+1)$-framed unknot by a simple $Y$-surgery (see Definition 
\ref{def:simpleY}). As a consequence, the Poincar\'e
sphere and the sphere $\mathbf{S}^{3}$, equipped with their unique 
spin structures, are $Y^{s}$-equivalent. Since $Y^{s}$-equivalence
is compatible with connected sums, \emph{we can assume} that $B=B'$. 

A theorem of Murakami and Nakanishi (\cite[Theorem 1.1]{MN}\footnote{In fact, 
the first reference is  Matveev,
but the proof in \cite{Mat} is not detailed.})
states that two ordered oriented 
links have identical linking matrices if and only if they are
$\Delta$-equivalent.
A $\Delta$-move is a certain unknotting operation, 
which is equivalent to surgery along a simple $Y$-graph.\\
Finally, from Corollary \ref{cor:simpleY},
we see that a simple $Y^{s}$-surgery between even surgery presentations
leaves the trivial characteristic solution fixed. 
We conclude that $(M,\sigma)$ is $Y^s$-equivalent to $(M',\sigma')$,
which completes the proof.  
\section{Applications}
\label{sec:application}
According to Theorem \ref{th:main}, 
two connected closed  spin $3$-manifolds are 
$Y^{s}$-equivalent if and only if they are $Y$-equivalent as \emph{plain} 
$3$-manifolds and their Rochlin invariants are identical modulo $8$.
In other words, while studying the degree $0$ part of Goussarov-Habiro theory,
the spin problem can be ``factored out''.

Now, given a closed connected oriented $3$-manifold, one can wonder whether  
all of its spin structures are $Y^{s}$-equivalent one to another. This has been
verified to be true in the case of 
$\mathbf{S}^{1} \times \mathbf{S}^{1} \times
\mathbf{S}^{1}$ by a direct calculation (Example \ref{ex:torus}). More
generally we have:
\begin{corollary}
\label{cor:corofmainth}
Let $M$ be a  connected oriented closed $3$-manifold such that 
$H_{1}(M)$ has no cyclic direct summand of order $2$ or $4$.
Then, all spin structures of $M$ are $Y^{s}$-equivalent one to another.
\end{corollary}
\begin{proof} 
This follows directly from Theorem \ref{th:main} and 
Corollary \ref{cor:top_noloworder}.
\end{proof}
On the contrary, we have:
\begin{example}
The two spin structures of
$\mathbf{RP}^{3}$ are not $Y^{s}$-equivalent, for the Rochlin function of
$\mathbf{RP}^{3}$ takes $1$ and $-1$ as values.
\end{example}
\bibliographystyle{amsalpha}

\begin{thebibliography}{99}
\newcommand{\no}{$\textrm{n}^{\circ}$}
\bibitem[Bl]{Bl}
   C. Blanchet,
   \emph{Invariants on three-manifolds with spin-structure},
   Comment. Math. Helvitici $\mathbf{67}$ (1992), 406--427.
\bibitem[CM]{CM}
   T.D. Cochran, P. Melvin,
   \emph{Finite type invariants of $3$-manifolds}, 
   Invent. Math. $\mathbf{140}$ (2000), 45--100.
\bibitem[De1]{DelI}
   F. Deloup,
   \emph{Linking forms, reciprocity for Gauss sums and 
   invariants of $3$-manifolds},
   Trans. of the A.M.S. $\mathbf{351}$ \no 5 (1999), 1895--1918.
\bibitem[De2]{DelII}
   \bysame,
   \emph{Reciprocity for Gauss sums and invariants of $3$-manifolds},
   Th\`ese de doctorat (1997), Universit\'e de Strasbourg.
\bibitem[Du]{Dur}
   A.H. Durfee,
   \emph{Bilinear and quadratic forms on torsion modules},
   Adv. in Math. $\mathbf{25}$ (1977), 133--164.  
\bibitem[GGP]{GGP}
   S. Garoufalidis, M. Goussarov, M. Polyak,
   \emph{Calculus of clovers and FTI of $3$-manifolds}, 
   Geometry and Topology $\mathbf{5}$ (2001), 75--108.
\bibitem[Gi]{Gil}
   C. Gille, 
   \emph{Sur certains invariants r\'ecents en topologie
   de dimension 3},
   Th\`ese de Doctorat (1998), Universit\'e de Nantes.
\bibitem[Go]{Gou}
   M. Goussarov,
   \emph{Finite type invariants and n-equivalence of $3$-manifolds},
   Compt. Rend. Ac. Sc. Paris $\mathbf{329}$ S\'erie I (1999), 517--522.
\bibitem[Ha]{Hab}
   K. Habiro,
   \emph{Claspers and finite type invariants of links},
   Geometry and Topology 
   $\mathbf{4}$ (2000), 1--83.
\bibitem[Jo]{Joh}
   D. Johnson,
   \emph{Spin structures and quadratic forms on surfaces},
   J. London Math. Soc. (2) $\mathbf{22}$ (1980), 365--373.
\bibitem[Ka]{Kap}
   S.J. Kaplan, 
   \emph{Constructing framed 4-manifolds with given almost
   framed boundaries}, 
   Trans. of the A.M.S. $\mathbf{254}$ (1979), 237--263.
\bibitem[KK]{KK}
   A. Kawauchi, S. Kojima,
   \emph{Algebraic Classification of Linking Pairings},
   Math. Ann. $\mathbf{253}$ (1980), 29--42.
\bibitem[Ki]{Kir}
   R.C. Kirby,
   \emph{The topology of $4$-manifolds},
   LNM $\mathbf{1374}$, Springer-Verlag (1991).
\bibitem[LL]{LL}
   J. Lannes, F. Latour,
   \emph{Signature modulo $8$ des vari\'et\'es de dimension $4k$ dont le bord
   est stablement parall\'elis\'e},
   Compt. Rend. Ac. Sc. Paris $\mathbf{279}$ S\'erie A (1974), 705--707.
\bibitem[Li]{Lic}
   W.B.R. Lickorish,
   \emph{A representation of orientable combinatorial $3$-manifolds},
   Ann. of Math. $\mathbf{76}$ (1962), 531--540.
\bibitem[Ma]{Mat}
   S.V. Matveev,
   \emph{Generalized Surgery of three-dimensional manifolds and representations
   of homology spheres}, Mat. Zametki $\mathbf{42}$ \no 2 (1987), 268--278
   (English translation in: Math. Notices Acad. Sci. USSR, $\mathbf{42:2}$).
\bibitem[MH]{MH}
   J. Milnor, D. Husemoller,
   \emph{Symmetric bilinear forms},
   Ergebnisse der Math. 73, Berlin, Heidelberg, New York (1973).
\bibitem[MS]{MS}
   J. Morgan, D. Sullivan,
   \emph{The transversality characteristic class and linking cycles in
   surgery theory},
   Ann. of Math. II Ser. 99 (1974), 463--544.
\bibitem[MN]{MN} 
   H. Murakami, Y. Nakanishi,
   \emph{On a certain move generating link-homology},
   Math. Ann. $\mathbf{284}$ (1989), 75--89.
\bibitem[Tu]{Tur}
   V.G. Turaev,
    \emph{Cohomology rings, linking forms and invariants of spin structure of
   three-dimensional manifolds}, Math. USSR Sbornik 
   $\mathbf{48}$ \no 1(1984), 65--79.
\bibitem[VdB]{VDB}
   F. Van der Blij,
   \emph{An invariant of quadratic forms modulo 8},
   Indag. Math. $\mathbf{21}$ (1959), 291--293.
\bibitem[Wa1]{WalI}
   C.T.C. Wall,
   \emph{Quadratic forms on finite groups, and related topics},
   Topology $\mathbf{2}$ (1964), 281--298.
\bibitem[Wa2]{WalII}
   \bysame,
   \emph{Quadratic forms on finite groups II},
   Bull. London Math. Soc. $\mathbf{4}$ (1972), 156--160.    
\bibitem[Wa3]{WalIII}
   \bysame,
   \emph{Non-additivity of the signature},
   Invent. Math. $\mathbf{7}$ (1969), 269--274. 
\end{thebibliography}

\vspace{0.5cm}
\noindent
\footnotesize{Commutative diagrams were drawn with Paul Taylor's package.}
\vspace{2cm}
\end{document}